\documentclass[11pt,reqno]{amsart}

\usepackage[text={160mm,240mm},centering]{geometry}            
\usepackage{tikz}
\usepackage{tikz-3dplot}
\usetikzlibrary{arrows.meta}
\usetikzlibrary{commutative-diagrams}
\usetikzlibrary{cd}
\usepackage{amssymb}
\usepackage{amsmath}
\usepackage{longtable}
\usepackage{textcomp}
\newtheorem{theorem}{Theorem}[section]
\newtheorem{lemma}[theorem]{Lemma}
\newtheorem{proposition}[theorem]{Proposition}
\newtheorem{conjecture}[theorem]{Conjecture}
\newtheorem{corollary}[theorem]{Corollary}

\theoremstyle{definition}
\newtheorem{definition}[theorem]{Definition}
\newtheorem{example}[theorem]{Example}

\theoremstyle{remark}
\newtheorem{remark}[theorem]{Remark}

\numberwithin{equation}{section}

\begin{document}

\title[unimodal hypersurface singularities]{Extended Tjurina number jumps and unimodal hypersurface singularities in positive characteristic}

\author{Hongrui Ma}
\address{Department of Mathematical Sciences,
Tsinghua University, Beijing, 100084, P. R. China.}
\email{mhr23@mails.tsinghua.edu.cn}

\author{Aoyu Ying}
\address{Zhili College, Tsinghua University, Beijing, 100084,
P. R. China.}
\email{yay23@mails.tsinghua.edu.cn }

\author{Huaiqing Zuo}
\address{Department of Mathematical Sciences,
Tsinghua University,
Beijing, 100084, P. R. China.}
\email{hqzuo@mail.tsinghua.edu.cn}





\begin{abstract}
This paper generalizes existing methods to derive stronger bounds on the modality of hypersurface singularities. Our results demonstrate that each sudden jump in the extended Tjurina number necessarily increases the modality. Furthermore, we provide a full classification of unimodal isolated hypersurface singularities in characteristic $p>3$ under contact equivalence.

	
\end{abstract}

\maketitle

\section{Introduction}
Throughout this paper, let $K$ be an algebraically closed field of arbitrary characteristic. We denote by $R=K[[\mathbf{x}]]=K[[x_1,\dots,x_n]]$  the formal power series ring and by $\mathfrak{m}=\langle x_1,\dots,x_n \rangle$ its maximal ideal. For an isolated hypersurface singularity we mean a power series $f \in R$ for which the Tjurina number $\tau(f)=\dim_K R/\langle f,\frac{\partial f}{\partial x_1},\dots, \frac{\partial f}{\partial x_n} \rangle$ is finite. We denote the extended Tjurina number $\tau^e(f)=\dim_K R/\langle f \rangle+\mathfrak{m}\cdot\langle\frac{\partial f}{\partial x_1},\dots, \frac{\partial f}{\partial x_n} \rangle$. Then $\tau(f)$ is finite if and only if $\tau^e(f)$ is finite.

The classification of singularities represents a fundamental challenge and a primary aim of singularity theory. In the classification of singularities, there are two equivalence relations: contact equivalence and right equivalence. Two power series $f,g \in R$ are contact equivalent if there exists a unit $U \in R^{\times} $ and an automorphism $ \phi \in Aut(R)$ such that $g=U\cdot \phi(f)$.

The modality of singularities for real and complex hypersurfaces was first introduced by Arnold in \cite{arnold}: the modality of a point $x \in X$ under the action of a Lie group $G$ on a manifold $X$ is the smallest $m$ such that a sufficiently small neighborhood of $x$ may be covered by a finite number of orbit families of $m$ parameters. Arnold \cite{arnold-class-C} completed the classification of hypersurface singularities with small modality over $\mathbb{C}$ under right equivalence. Subsequently, Wall \cite{Wall} established the classification of unimodal hypersurface singularities under contact equivalence.

In their work \cite{1990_simple}, Greuel and Kroning classified hypersurface singularities of finite deformation type (modality 0) over fields of positive characteristic under contact equivalence, employing finite determinacy theory. Later, Boubakri, Greuel, and Markwig \cite{2010_finite_deter} refined the finite determinacy theorem in 2010, which has since become a fundamental tool in classification problems.

In 2016, Greuel and Nguyen \cite{right-simple} extended the concept of modality of hypersurface singularities to arbitrary algebraically closed fields by developing an algebraic formulation. Their work established a fundamental theorem providing explicit bounds for modality in this generalized setting:

\begin{theorem}\label{1}
    Assume $X$ is irreducible,  for every $x \in X$, let $G\text{-modality}(x)$ be the modality of $x$  under $G$ (see Definition \ref{def-mod}). Then
    $$G\text{-modality}(x) \geq \dim X-\dim G.$$
\end{theorem}
Greuel and Nguyen further classified hypersurface singularities 
of modality 0 in positive characteristic under right equivalence  \cite{right-simple}. Subsequently, Nguyen \cite{right_unimodal} 
extended this classification to singularities of modality 1 and 2 under right equivalence. Recently, Pham, Pfister, and Greuel extended the concept of modality to isolated complete intersection singularities (ICIS) over arbitrary algebraically closed fields \cite{icissimple}. Building on this work, they classified ICIS of modality 0 under contact equivalence in positive characteristic. Shortly thereafter, Ma, Yau, and Zuo provided a classification for ICIS of modality 1 under the same equivalence in positive characteristic \cite{ma2025}.

To apply Theorem~\ref{1}, one must choose a suitable jet space $X$ for a given 
power series jet, typically requiring significant computation. This process 
becomes tractable under right equivalence owing to the finite classes
in characteristic~$p$. The case of contact equivalence differs markedly: 
the germ $x^2 + y^k$ has the right modality 0 only when $k \leq p-1$, 
but it has the contact modality 0 for all $k \geq 2$.

In this paper, we generalize Theorem~\ref{1} to Theorem \ref{better-mod-boundary} and obtain sharper bounds on modality. Under some further assumptions, we obtain Theorem \ref{main-thm}, which establishes a fundamental connection between the modality of a family of hypersurface singularities and their extended Tjurina numbers.

A key observation is that sudden jumps in the extended Tjurina number may occur for families of singularities over fields of positive characteristic. For example, consider the family $f_k = x^3 + xy^{k} + y^{34}$ ($23\leq k \leq 32$) over $\mathbb{C}$, where $\tau^e(f_k) = k + 35$ grows linearly. However, in the field of characteristic 5, we find:
\[
\tau^e(f_{26}) = 68 \neq 61 \quad \text{(unexpected jump)}
\]
while $\tau^e(f_{25}) = 60$ and $\tau^e(f_{27}) = 62$ remain consistent with the complex case. This phenomenon - where certain singularities exhibit extended Tjurina numbers strictly greater than their neighbors - is what we call a \emph{sudden jump} of the extended Tjurina number. Our results demonstrate that each such jump necessarily increases the modality.

Building on this framework, we complete the classification of unimodal hypersurface singularities  in characteristic $p > 3$ under contact equivalence. The classification, presented in Theorem~\ref{final-class}, is more intricate than in the complex case due to the need to account for these sudden jumps of extended Tjurina numbers. For the cases of small characteristic, however, things become more complicated since the orbit map $o: G \rightarrow G\cdot f$ may not be separable, and there are still lots of works to be done.

{\bf Acknowledgement.}  We are grateful to the anonymous referees for their comments, which have helped improve this manuscript.  Zuo is supported by NSFC Grant 12271280 and BJNSF Grant 1252009.

\section{Contact equivalence and modality}
To fix notation, we first recall some key definitions.
\begin{definition}
    For a power series $f \in \mathfrak{m} \subset R$ we denote $tj(f)=\langle f,\frac{\partial f}{\partial x_1},\dots, \frac{\partial f}{\partial x_n} \rangle$ the Tjurina ideal of $f$. We call the associated algebra $T_f=R/tj(f)$ the Tjurina algebra. We call $f$ an isolated hypersurface singularity if the Tjurina number $\tau(f)$=$\dim_K T_f < \infty$.
\end{definition}

\begin{definition}
    The contact group $\mathcal{K}$ is defined as $$\mathcal{\mathcal{K}}=R^{\times} \rtimes Aut(R),$$ and the action of $\mathcal{\mathcal{K}}$ acting on $R$ is defined as $$(U,\phi,f) \mapsto U \cdot \phi (f),$$ with $U \in R^{\times},\ \phi \in Aut(R),\ f \in R$ and $$\phi(f)=f(\phi(\mathbf{x})),$$ where $\phi(\mathbf{x})=(\phi(\mathbf{x_1}),\dots,\phi(\mathbf{x_n}))$.
\end{definition}

Two isolated hypersurface singularities $f$ and $g \in K[[\mathbf{x}]]$ are called contact equivalent, denoted $f \sim_c g$ (or simply denoted $f \sim g$), if $g \in \mathcal{K}f$. 

Arnold introduced the definition of modality (see \cite{arnold}) over real or complex manifolds  as follows: The modality of a point $x \in X$ under the action of a Lie group $G$ on a manifold $X$ is the smallest $m$ such that a sufficiently small neighborhood of $x$
may be covered by a finite number of orbit families of $m$ parameters.

Greuel and Nguyen generalized the notion in the case of
hypersurface singularities over an algebraically closed field of arbitrary characteristic and gave a detailed discussion in \cite{right-simple}, \cite{phdclassification}. We collect some definitions here.

\begin{definition}\label{def-mod}
    Let $U\subset X$ be an open neighborhood of $x\in X$ and let $W$ be constructible in $X$. We introduce
\begin{eqnarray*}
\dim_x W&:=&\max\{\dim Z\ |\ Z \text{ is an irreducible component of }W \text{ containing }x \},\\
U(i)&:=&U_G(i):=\{ y \in U\ |\ \dim_y (U\cap G\cdot y ) =i\}, i\geq 0,\\
G\text{-mod}(U)&:=&\max_{i\geq 0}\{\dim U(i)-i\}.
\end{eqnarray*}
We define
$$G\text{-mod}(x):=\min \{G\text{-mod}(U)\ |\ U \text{ a neighborhood of x}\}$$
the modality of $x$ (in $X$) under $G$.
\end{definition}

For a function germ $f \in R^m$, denote by $J_k=R^m/\mathfrak{m}^{k+1}R^m$ the $k$-jet space of $R^m$. The $k$-jet of $f$ is the image in $J_k$, denoted by $j_k(f)$. Denote $\mathcal{K}_k=\{(j_k(U),j_k(\phi)) \mid U \in R^{\times},\ \phi \in Aut(R)\}$ as the $k$-jet contact group. Then the modality of $f$ under $\mathcal{K}$ is defined as the modality of a sufficiently large jet, denoted by $\mathcal{K}\text{-mod}(f)$.

Next we use the following facts from \cite{phdclassification} to give a criterion for non-unimodal.

\begin{proposition}\label{modality}Let an algebraic group $G$ act on a variety $X$.

    (1)If the subvariety $X' \subset X$ is invariant under $G$ and $x \in X'$, then $$G\text{-mod}(x) \mathrm{\ in}\ X \geq G\text{-mod}(x) \mathrm{\ in}\ X'.$$

    (2)Let additionally an algebraic group $G'$ act on a variety $X'$ and let $p:X \rightarrow X'$ be a morphism of varieties. $p$ is open and $$ G\cdot x \subset p^{-1}(G' \cdot p(x)),\ \forall x \in X.$$ Then $$G\text{-mod}(x) \geq G'\text{-mod}(p(x)),\ \forall x \in X.$$
    
    (3)If $X$ is irreducible, for $x \in X$, we have $$G\text{-mod}(x) \geq \mathrm{dim}X-\mathrm{dim}G.$$
\end{proposition}

\begin{proposition}\label{nl-bound}
    Let $f \in K[[x_1,\dots,x_n]]$ be a unimodal (i.e. of modality 1) isolated hypersurface singularity. Let $\mathrm{ord}(f)=l$. Then one of the following holds:\\
    (i) $n \geq 4,\ l=2$;\\
    (ii) $n=3,\ l \leq 3$;\\
    (iii) $n=2,\ l \leq 4$.
\end{proposition}
\begin{proof}
   Choose $k$ sufficiently large and let $X=\mathfrak{m}^l/\mathfrak{m}^{k+1}$. It follows from Proposition \ref{modality}(1) that $$1=\mathcal{K}\text{-mod}(f)=\mathcal{K}_k-mod (f) \mathrm{\ in\ } J_k \geq \mathcal{K}_k-mod (f) \mathrm{\ in\ } X.$$
    Let $X'=\mathfrak{m}^l/\mathfrak{m}^{l+1}$. The action of $\mathcal{K}_k$ on $X$ induces the action of the algebraic group $\mathcal{K}'=I\times GL(n,K)$ on $X'$, and it can easily be checked that $p:X \rightarrow X'$ is open and $\mathcal{K}_k \cdot f \subset p^{-1}(\mathcal{K}' \cdot p(f))$. Then by Proposition \ref{modality}(2) we have $$\mathcal{K}_k-mod (f) \mathrm{\ in\ } X \geq \mathcal{K}'-mod (p(f)) \mathrm{\ in\ } X' .$$ Therefore by Proposition \ref{modality}(3) we have $$1 \geq \mathcal{K}'-mod (p(f)) \mathrm{\ in\ } X' \geq  \mathrm{dim}X'-\mathrm{dim}\mathcal{K}'. $$
    
    Calculation shows that $$\mathrm{dim}X'=\binom{n-1+l}{l} \mathrm{\ and\ }\mathrm{dim}\mathcal{K}'=n^2. $$

    Thus $1 \geq \binom{n-1+l}{l}-n^2$. The solution is what we want.
\end{proof}

\section{Classification methods}\label{method}
The finite determinacy theorem plays a crucial role in the proof of the classification theorem.
\begin{theorem}[\cite{finitedeter}]\label{finite-deter}
    Let $f \in \mathfrak{m}^2$. If there exists a natural number $k \in \mathbb{N}$ such that $$\mathfrak{m}^{k+2} \subset \mathfrak{m} \cdot \widetilde{T}_f(\mathcal{K}f),$$ then $f$ is $(2k-\mathrm{ord}(f)+2)$-determined, where $$\widetilde{T}_f(\mathcal{K}f)=\langle f \rangle+\mathfrak{m}\cdot \langle \frac{\partial f}{\partial x_1}, \dots, \frac{\partial f}{\partial x_n}\rangle $$ is the tangent image. That is, for any $g \in R^m$ with $j_{2k-\mathrm{ord}(f)+2}(g)=j_{2k-\mathrm{ord}(f)+2}(f)$, we always have $g \sim f$.
\end{theorem}

\begin{remark}
    We denote $$T^e_f=R/\widetilde{T}_f(\mathcal{K}f)=K[[\mathbf{x}]]/(\langle f \rangle+\mathfrak{m}\cdot \langle \frac{\partial f}{\partial x_1}, \dots, \frac{\partial f}{\partial x_n}\rangle)$$ the expanded Tjurina algebra $\tau^e(f)=\dim_K T^e_f $ the extended Tjurina number, which will be mentioned in the following sections. Note that $\tau(f) < \infty$ if and only if $\tau^e(f) < \infty.$
\end{remark}

The following method is the generalization of the finite determinacy theorem, which is developed in \cite{normal}. We collect the main results here.

Given $\mathbb{Q}$-linear independent weight vectors $w_i\in\mathbb{Q}_{>0}^n$ with positive entries,
   $i=1,\ldots,k$, they define linear functions
   \begin{displaymath}
     \lambda_i:\mathbb{R}^n\longrightarrow\mathbb{R}:r\mapsto w_i\cdot r:=\sum_{j=1}^n
     w_{i,j}\cdot r_j,
   \end{displaymath}
   which induces
   \begin{displaymath}
     \lambda:\mathbb{R}^n\longrightarrow\mathbb{R}:r\mapsto\min\{\lambda_1(r),\ldots,\lambda_k(r)\}.
   \end{displaymath}
   The set
   \begin{displaymath}
     P_\lambda=\{r\in\mathbb{R}_{\geq 0}^n\;|\;\lambda(r)=1\}
   \end{displaymath}
   is a compact rational polytope of dimension $n-1$ in the positive orthant $\mathbb{R}_{\geq 0}^n$,
   and its facets are given by
   \begin{displaymath}
     \Delta_i=\{r\in P_\lambda\;|\;\lambda_i(r)=1\}.
   \end{displaymath}
   Such sets are called \emph{$C$-polytopes}. Thus, $\mathbb{Q}$-linear independent weight vectors define $C$-polytopes. Conversely, given a $C$-polytope $P$, we can get a set of $\mathbb{Q}$-linear independent weight vectors.

   For a $C$-polytope $P$, we denote $N_P$ the lowest common multiple
   of the denominators of all entries in the weight vectors
   corresponding to $P$. Then we can define a valuation on $K[[\mathbf{x}]]$ by 
   $$v_P(f):=\mathrm{min}_{\alpha}\{N_p \cdot \lambda_P(\alpha) \mid \alpha \in \mathrm{supp}(f)\}$$ for a power series $f=\sum_\alpha a_\alpha \mathbf{x}^{\alpha}\in K[[\mathbf{x}]]$, where $\mathrm{supp}(f)=\{\alpha \in \mathbb{N}^n \mid a_\alpha \neq 0\}$. Suppose that the corresponding weight vectors of $P$ are $w_i,\ i=1,\dots, k$, we define $$v_i(f):=\mathrm{min}\{N_P \cdot \lambda_i(\alpha) \mid \alpha \in \mathrm{supp}(f)\}. $$ Then $v_P$ satisfies $$v_P(f \cdot g) \geq v_P(f)+v_P(g),\ v_P(f+g) \geq \mathrm{min}\{v_P(f),v_P(g)\}$$ and 
   \begin{equation}
     \label{eq:vp=vd}
     v_P(f\cdot g)=v_P(f)+v_P(g)
     \;\;\;\Longleftrightarrow\;\;\;
     v_P(f)=v_i(f)\;\;\mbox{ and }\;\;
     v_P(g)=v_i(g)     
   \end{equation} for some $i$.

   Note that for a power series $f \in K[[\mathbf{x}]]$ as above, the Newton diagram $\Gamma(f)$ of $f$ is a $C$-polytope if and only if $f$ is a \emph{convenient} power series, i.e. if the support of $f$ contains a point on each coordinate axis. We denote $v_{\Gamma(f)}$ simply as $v_f$. In this case, we have $v_f(f)=v_i(f)$ for all $i=1,\dots,k.$ If $f$ is not convenient, we usually expand the Newton diagram in a suitable way to obtain the $C$-polytope $P$. 

\begin{example}
    Let $f=x^3+xy^r+y^s+yz^2+z^3 \in K[[x,y,z]]$. Then $f$ is convenient. The Newton diagram $\Gamma(f)$ is shown in Figure \ref{fig-newton} of case $r=3,s=5$.
    \\
    \begin{figure}[htbp]
    \centering
    \begin{tikzpicture}[scale=2]

\coordinate (O) at (0,0,0);
\coordinate (A) at (0,0,1.5);
\coordinate (B) at (2.5,0,0);
\coordinate (C) at (0,1.5,0);
\coordinate (D) at (1.5,0,0.5);
\coordinate (E) at (0.5,1,0);

\draw[thick,->] (0,0,0) -- (3,0,0) node[anchor=north east]{$y$};
\draw[thick,->] (0,0,0) -- (0,2,0) node[anchor=north west]{$z$};
\draw[thick,->] (0,0,0) -- (0,0,2) node[anchor=south]{$x$};

\draw[thick] (A) -- (D);
\draw[thick] (D) -- (B);
\draw[thick] (A) -- (C);
\draw[thick] (C) -- (E);
\draw[thick] (B) -- (E);
\draw[dashed] (A) -- (E);
\draw[dashed] (D) -- (E);

\foreach \point in {O, A, B, C, D, E}
    \fill[black] (\point) circle (1.5pt);
\end{tikzpicture}
\caption{The Newton diagram of $x^3+xy^3+y^5+yz^2+z^3$.}
\label{fig-newton}
\end{figure}

The corresponding weight vectors are $$w_1=(2rs,4s,3rs-2s),\ w_2=(6rs-6r^2,6r,3rs-3r),\ w_3=(2rs,2rs,2rs)$$ and $v_f(f)=v_i(f)=6rs$ for $i=1,2,3$.
\end{example}

\begin{example}
    Let $f=x^3+xy^6 \in K[[x,y]]$. Then $f$ is not convenient. We choose $P$ expanded from the Newton diagram given by $(0,3),(6,1),(9,0)$ (the expanding point) as shown in Figure \ref{fig-expanding-newton}.
        \\
    \begin{figure}[htbp]
    \centering
    \begin{tikzpicture}[scale=2]

\coordinate (O) at (0,0);
\coordinate (A) at (0,1);
\coordinate (B) at (2,0.33);
\coordinate (C) at (3,0);

\draw[thick,->] (0,0) -- (3.5,0) node[anchor=north east]{$x$};
\draw[thick,->] (0,0) -- (0,1.5) node[anchor=north west]{$y$};

\draw[thick] (A) -- (B);
\draw[dashed] (B) -- (C);

\foreach \point in {O, A, B, C}
    \fill[black] (\point) circle (1.5pt);
\end{tikzpicture}
\caption{$C$-polytope expanded from the Newton diagram of $x^3+xy^6$.}
\label{fig-expanding-newton}
\end{figure}

\end{example}

We can extend $v_P$ to $\mathrm{Der}_K(K[[\mathbf{x}]])$ as following: for $$\xi=\sum_{i=1}^n\sum_{\alpha\in\mathbb{N}^n} a_{i,\alpha}\cdot
     \mathbf{x}^{\alpha} \cdot \partial_{x_i} \in \mathrm{Der}_K(K[[\mathbf{x}]]),$$ let $$v_P(\xi)=\mathrm{min}\{\lambda_P(\alpha-e_i) \mid a_{i,\alpha} \neq 0\}.$$ It follows that $$v_P(\xi f) \geq v_P(\xi)+v_P(f).$$

For a $C$-polytope $P$, taking the filtration induced by $v_P$, denoted by $F_d$, i.e. $F_d=\{h \in K[[\mathbf{x}]] \mid v_P(h) \geq d\}$. 

Furthermore, for $f$ is a hypersurface singularity, define
$$tj(f)_d:=\{h=g\cdot f+\xi f\;|\;g\in K[[\mathbf{x}]],\xi\in\mathrm{Der}_K(K[[\mathbf{x}]]),v_P(h) \geq d\}$$ the graded Tjurina ideal and

$$tj^{AC}(f)_d:=\{h=g\cdot f+\xi f\;|\;\min\{v_P(g)+v_P(f),v_P(\xi)+v_P(f)\}\geq d\}$$ the AC-graded Tjurina ideal.

Then we have the graded algebras
$$gr_P(T_f):=\bigoplus_{d\geq 0} F_d\big/\big(tj(f)_d+F_{d+1}\big) \cong K[[\mathbf{x}]]/tj(f)= T_f$$ and 
$$gr_P^{AC}(T_f):=\bigoplus_{d\geq 0} F_d\big/\big(tj^{AC}(f)_d+F_{d+1}\big).$$ Clearly we have $$gr_P^{AC}(T_f) \twoheadrightarrow T_f.$$

\begin{definition}
    A monomial basis of $gr_P^{AC}(T_f)$ 
     is called a regular basis for $T_f$.
\end{definition}

Given any $C$-polytope $P$ and a power series $f\in K[[\mathbf{x}]]$, we call
$$in_P(f)=\sum_{\stackrel{\lambda_P(\alpha) \text{ minimal}}{\alpha\in Supp(f)}}a_\alpha {\bf x}^\alpha$$ the initial part of $f$. One can show $gr_P^{AC}(T_f)=gr_P^{AC}(T_{in_P(f)})$ (\cite{normal} Lemma 3.8). Same as above, we write $in_f(f)$ instead of $in_{\Gamma(f)}(f)$ when we choose $P=\Gamma(f)$.

\begin{theorem}[\cite{normal} Theorem 4.5]\label{normal-thm}
    Let $f\in\mathfrak{m}$, $P$ be a $C$-polytope and
     $B=\{x^\alpha\;|\;\alpha\in\Lambda\}$ a regular basis for
     $T_{in_P(f)}$. 
     If 
     \begin{equation}\label{finite-deter-condition}
         \mathfrak{m}^{k+2} \cdot R^m \subset \mathfrak{m} \cdot \widetilde{T}_f(\mathcal{K}f),
     \end{equation}
      then 
     \begin{equation}\label{eq:nfce:1}
       f\;\sim\;in_P(f)+\sum_{\alpha\in\Lambda_f} c_\alpha {\bf x}^\alpha.
     \end{equation}
     for suitable $c_\alpha\in K$, where $\Lambda_f$ is the finite set
     $$
       \Lambda_f=\big\{\alpha\in\Lambda\;\big|\;\deg({\bf x}^\alpha)\leq
       2k-\mathrm{ord}(f)+2,\; v_P({\bf x}^\alpha)> v_P\big(in_P(f)\big)\big\}.$$
\end{theorem}

However, we can hardly find a suitable $k$ satisfying \ref{finite-deter-condition} if we don't know the normal form of $f$. We have the following corollary avoiding condition \ref{finite-deter-condition}.
\begin{corollary}[\cite{normal} Corollary 4.7] \label{normal-cor}
     Let $P$ be a $C$-polytope and $f\in\mathfrak{m}$ be a power series such
     that $in_P(f)$ satisfies $\mathrm{dim}gr_P^{AC}(T_{in_P(f)})<\infty$,  then $f$ is finitely determined, and
     \begin{displaymath}
       f\;\sim\;in_P(f)+\sum_{\stackrel{{\bf x}^\alpha\in
           B}{v_P({\bf x}^\alpha)>d}} c_\alpha {\bf x}^\alpha
     \end{displaymath}
     for suitable $c_\alpha\in K$, where $B$ is a finite regular
     basis for $T_{in_P(f)}$ and $d=v_P(in_P(f))$. 
   \end{corollary}

To avoid redundant coefficients $c_\alpha$, we often employ the following implicit function theorem.

\begin{theorem}[\cite{singular-introduction} Theorem 6.2.17]\label{implicit-function}
    Let $\mathcal{K}$ be a field and $F \in \mathcal{K}[[x_1,\dots,x_n,y]]$ such that
    \begin{equation}
        F(x_1,\dots,x_n,0) \in \langle x_1,\dots,x_n \rangle,\ \frac{\partial F}{\partial y}(x_1,\dots,x_n,0) \notin \langle x_1,\dots,x_n \rangle,
    \end{equation}
    then there exists a unique $y(x_1,\dots,x_n) \in \langle x_1,\dots,x_n \rangle \mathcal{K}[[x_1,\dots,x_n]]$ such that $$F(x_1,\dots,x_n,y(x_1,\dots,x_n))=0.$$
\end{theorem}
We will show the use of Theorem \ref{implicit-function} in Section \ref{classification p}.

In fact, Theorem \ref{normal-thm} and \ref{normal-cor} give us a better bound of finite determinacy than Theorem \ref{finite-deter}. 
\begin{corollary}[\cite{normal} Corollary 4.9]\label{normal-deter}
    Let $P$ be a $C$-polytope and $f\in\mathfrak{m}$ be a power series such
     that $in_P(f)$ satisfies $\mathrm{dim}gr_P^{AC}(T_{in_P(f)})<\infty$. Let $B$ be a regular basis of $T_{in_P(f)}$. Then $$d:=\max_{\mathbf{x}^\alpha \in B} \{v_p(in_p(f)),v_p(\mathbf{x}^\alpha)\} $$ is finite and $f \sim g$ for every $g \in R$ with $v_p(f-g)>d.$ Moreover, if $\mathfrak{m}^{k+1} \in F_{d+1}$, then $f$ is $k$-determined.
\end{corollary}

For a given $k$-jet of $f$, the following theorem from \cite{complete-unimodal-plane} can be used to confirm $in_P(f)$.  In \cite{ma2025}, the authors have modified some notation to match the case of positive characteristic fields.

\begin{theorem}\label{complete-cor}  
    Let $f \in J_k$ be a $k$-jet of weighted homogeneous type w.r.t. $(a_1,\dots,a_n;d)$. That is, $f$ satisfies $$f(t^{a_1}x_1,\dots,t^{a_n}x_n)=t^{d}f(x_1,\dots,x_n).$$ Moreover, assume  
    \begin{equation}\label{condition-1.2}
        d<(k+1)\mathrm{min}(a_j) \mathrm{\ or\ }d>(k+1)\mathrm{max}(a_j).
    \end{equation}
    Denote $P_{k,l}=\mathfrak{m}^{k+1}/\mathfrak{m}^{l+1}$ as a linear space. Let $C \subset P_{k,l}$ be a linear subspace of $P_{k,l}$ satisfying $$P_{k,l} \subset C+\widetilde{T}_f(\mathcal{K}_{l}f)\cap P_{k,l},$$ we call $C$ a complete transversal. This complete transverse has the following property: every $g \in J_l$ of the same $k$-jet with $f$ is in the same $\mathcal{K}_{l}$-orbit as some $l$-jet of the form $f+c$, for some $c \in C$.
\end{theorem}

Same method can also be used to find modality. See also \cite{ma2025}.

Let $C$ be a complete transversal of $f$ in $J_l\ (l>k),$ for $a \in C$, we define
\begin{equation}
    \mathrm{cod}(f+a)=\mathrm{comdimension\ of\ }\widetilde{T}_f(\mathcal{K}_{l}f) \cap P_{k,l} \mathrm{\ in\ }P_{k,l}
\end{equation}
and
\begin{equation}
    \mathrm{cod}_0(f)=\mathrm{inf}_{a \in C}\{\mathrm{cod}(f+a)\}.
\end{equation}
Note that there exists a Zariski open subset $U \subset C$ such that $\mathrm{cod}(f+a)=\mathrm{cod}_0(f)$ if and only if $a \in U.$
\begin{theorem}\label{cod-modality}
    Let $f$ be defined as above. Then for $a \in U,$ $f+a$ has modality $\mathrm{cod}_0(f)$ in $J_l(f)$ under the action of the subgroup $\mathcal{K}_l(f)$ of $\mathcal{K}_l$ which stabilize $f$. In particular, any jet $h$ in $J_l(f)$ has $\mathcal{K}_l(f) \text{-mod}(h) \geq \mathrm{cod}_0(f)$ in $J_l$.
\end{theorem}

We will show the use of Theorem \ref{normal-thm} $\sim$ \ref{cod-modality} in Section \ref{classification p}.

\section{A new criterion of modality of hypersurface singularity}\label{section-new-mod}
By \cite{Rosenlicht} Theorem 2, for an algebraic group $G$ acting on a variety $X$, there exists an open dense set $X_1 \subset X$, which is invariant under $G$,  such that $X_1/G$ is a geometric quotient. In particular, $X_1/G$ is an algebraic variety. If $X$ is irreducible, then $X_1/G$ is irreducible.

As we have mentioned above, Nguyen has shown in \cite{phdclassification} that $$G\text{-mod}(x) \geq \mathrm{dim}X-\mathrm{dim}G.$$ Using Rosenlicht's theorem, we can more precisely show that (with a little change of the original proof)
\begin{theorem}\label{better-mod-boundary}
    Let the algebraic group $G$ act on a variety $X$. If $X$ is irreducible, there exists a Zariski open subset $X_1 \subset X$, such that $$G\text{-mod}(x) \geq \mathrm{dim}X-\mathrm{dim}G\cdot x$$ for any $x \in X_1$.
\end{theorem}
\begin{proof}
    Let $U$ be an open neighborhood of $x \in X$ such that $G\text{-mod}(x)=G\text{-mod}(U).$ By definition, $$G\text{-mod}(U)=\max_{i\geq 0}\{\dim U(i)-i\}.$$
    We claim that: $$G\text{-mod}(U)=\max_{i\geq 0}\{\dim U(\leq i)-i\},$$ where $U(\leq i)=\{ y \in U\ |\ \dim_y (U\cap G\cdot y ) \leq i\}.$

    Note that $U(\leq i)=\bigcup_{j \leq i}U(j)$. The inequality $$\max_{i\geq 0}\{\dim U(i)-i\} \leq \max_{i\geq 0}\{\dim U(\leq i)-i\}$$ follows easily from $U(i) \subset U(\leq i)$. For the other side, we choose $i_0$ such that $$\max_{i\geq 0}\{\dim U(\leq i)-i\}=\dim U(\leq i_0)-i_0,$$ then we have
    \begin{equation}
        \begin{aligned}
            \max_{i\geq 0}\{\dim U(\leq i)-i\}&=\dim U(\leq i_0)-i_0 \\
            & =\max_{i \leq i_0}\dim U(i)-i_0 \\
            &=\max_{i\leq i_0}\{\dim U(i)-i\}\\
            &\leq \max_{i\geq 0}\{\dim U(i)-i\}.
        \end{aligned}
    \end{equation}
    The claim has been proved.

    By Rosenlicht's theorem, there exists an open dense set $X_1 \subset X$ such that $p_1:X_1 \rightarrow X_1/G$ is a dominant morphism of irreducible varieties. For every $y \in X_1$, we choose $i_1=\dim G\cdot y$. Thus, the set $$U_1=\{z \in X_1| \dim p_1^{-1}(p_1(z))\leq i_1\}=\{z \in X_1| \dim G\cdot z \leq i_1\}$$ is open and nonempty in $X_1$ by Chevalley's theorem, hence open in $X$. Therefore, 
    \begin{equation}
        \begin{aligned}
            G\text{-mod}(U)&=\max_{i\geq 0}\{\dim U(\leq i)-i\} \\
            &\geq \dim U(\leq i_1)-i_1 \\
            &\geq \dim (U \cap U_1)-\dim G \cdot y.
        \end{aligned}
    \end{equation}
    Since $X$ is irreducible, $U \cap U_1$ is a non-empty open subset of $X$, hence $\dim (U \cap U_1)=\dim X$, and we get $$G\text{-mod}(x) \geq \mathrm{dim}X-\mathrm{dim}G\cdot y$$ for every $x \in X$ and $y \in X_1$.
\end{proof}

\begin{remark}
     A more precise choice of \( i_1 \) will yield a better bound, as we will show in the next section.
\end{remark}

Next we consider the dimension of the orbit $G\cdot x$ for $x \in X$. The orbit map $o:G \rightarrow G \cdot x$ induces the tangent map $d_1o: T_eG \rightarrow T_x(G\cdot x)$. If $G$ is smooth, then $G \cdot x$ is smooth (cf. \cite{Milne_alg_group} Proposition 9.7), thus $\dim G\cdot x=\dim T_x(G \cdot x)$.

We introduce the definition of separable morphism.
\begin{definition}
    (i)We call the field extension $K/k$ separably generated if there exists a finite transcendence basis $\{x_i\}$ such that $K/k(\{x_i\})$ is separable.\\
    (ii) Let $\phi:X \rightarrow Y$ be a dominant morphism of irreducible algebraic varieties over $k$. Then it induces $\phi^\#: k(Y) \rightarrow k(X)$. We call $\phi$ a separable morphism if the extension $k(X)/\phi^{\#}(k(Y))$ is separably generated.
\end{definition}

We have the following theorem.
\begin{theorem}[\cite{actions-algebraic-group} Theorem 3.1]
    Let $G$ be an affine algebraic group, $X$ an algebraic $G$-variety and $x \in X$. Then the orbit $G \cdot x$ of $x$ is a non-singular algebraic variety of $X$. Moreover, the following are equivalent.\\
    (i) The orbit map $o:G \rightarrow G \cdot x$ is a separable morphism.\\
    (ii) The tangent map $d_1o: T_eG \rightarrow T_x(G\cdot x)$ is surjective.
\end{theorem}
Whether $d_1o$ is surjective or not, we denote the image of $d_1o$ as $\widetilde{T}_x(Gx)$, which has the same meaning as $\widetilde{T}_f(\mathcal{K}f)$ appearing in Theorem \ref{finite-deter}.

Now we set $f$ to be an isolated hypersurface singularity and $P=\Gamma(f)$ is the Newton diagram (or expanding from the Newton diagram) of $f$. 

Set $d=v_P(f),X=F_0/F_{l+1}$, where $l$ is an integer greater than $d$. Set $G=\mathcal{K}_l$, where the action of $G$ on $X$ is induced from the action of $\mathcal{K}$ on $R$. Specifically, we denote the natural projection $\pi: R=F_0 \rightarrow F_0/F_{l+1}$. The action of $G$ on $X$ is given by: 
    \begin{equation}
        \begin{aligned}
            &G \times X &\longrightarrow&X\\
            &((U,\phi),h) &\mapsto &\pi(U \cdot \phi(h)).
        \end{aligned}
    \end{equation}

\begin{proposition}\label{tangent-image-calculation}
    The tangent image $\widetilde{T}_f(Gf)=(\widetilde{T}_f(\mathcal{K}f)+ F_{l+1})/F_{l+1}$. 
\end{proposition}

\begin{proof}
    The orbit map $o:G \rightarrow G\cdot f$ is given by $$(j_l(U),j_l(\phi)) \rightarrow \pi(U \cdot \phi(f)).$$

    Each element of $T_eG$ can be written as $(j_l(1+\epsilon U),j_l(id_R+\epsilon \phi))$, where $\epsilon^2=0$. We write $$\phi:(x_1,\dots,x_n) \mapsto (x_1+\phi_1,\dots,x_n+\phi_n).$$ Acting on $f$, we get $$\pi((1+\epsilon U)\cdot f(x_1+\phi_1,\dots,x_n+\phi_n)).$$

    Using the Taylor expansion, we have 
     \begin{equation}
        (1+\epsilon U)\cdot f(x_1+\phi_1,\dots,x_n+\phi_n)=f(\mathbf{x})+\epsilon Uf(\mathbf{x})+\epsilon \sum_i \frac{\partial f}{\partial x_i}\phi_i.
    \end{equation}
    Therefore, the image of the tangent map $d_1o$ is generated by the image of $Uf,\sum\frac{\partial h}{\partial x_i}\phi_i$ under $\pi$, which coincides with 
    $$\{Uf+\sum_i\phi_i\frac{\partial f}{\partial x_i}+F_{l+1}\}/F_{l+1}=(\widetilde{T}_f(\mathcal{K}f)+ F_{l+1})/F_{l+1}.$$
\end{proof}

\begin{corollary}\label{dim-Gf-cor}
    Denote a basis of $T^e_f$ by $B^e_f$ and choose $l=\max_{\mathbf{x}^\alpha\in B^e_f}\{v_P(\mathbf{x^\alpha})\}$. Let $X=F_0/F_{l+1},\ \mathcal{K}=G_l$ be as above. If the orbit map $o:G \rightarrow G \cdot f$ is separable, then $\dim G\cdot f=\dim \widetilde{T}_f(Gf)=\dim X- \tau^e(f)$.
\end{corollary}
\begin{proof}
    In this case, $\widetilde{T}_f(Gf) \subset F_{l+1}$ by the definition of $l$. Therefore, by Proposition \ref{tangent-image-calculation}, 
    $$\dim G\cdot f=\dim \widetilde{T}_f(Gf)=\dim \widetilde{T}_f(\mathcal{K}f)=\dim X-\dim X/\widetilde{T}_f(\mathcal{K}f)=\dim X-\tau^e(f). $$
\end{proof}

\begin{remark}\label{remark-separable}
    (i) If char$K=0$, then the orbit map is always separable since the field extension is always separable over a characteristic $0$ field. Hence, the result in Corollary \ref{dim-Gf-cor} always holds.\\
    (ii) If char$K=p>0$, then there exists $f$ such that some orbit maps may not be separable. See \cite{finitedeter} Example 2.9. However, each of the counterexamples given satisfies $p \mid ord  (f)$. In fact, we can show that for $f$ of the form $x^p+\mathfrak{m}^{p+1}$, the orbit map $o:\mathcal{K}_k \rightarrow \mathcal{K}_k \cdot f$ cannot be separable: write $\phi(x)=a_{11}x+a_{12}y+a_{21}x^2+\dots$, then $K(a_{11}) \subset K(\mathcal{K}_k)$ and $K(a_{11}^p) \subset K(\mathcal{K}_k\cdot f)$, then $K(a_{11})/K(a_{11}^p)$ is not separable. But if we choose $f$ such that $ord(f) \leq 4$ in the field of characteristic $p$ greater than $5$, and set the space $X=F_0/F_{l+1}$, every example we calculate shows that the orbit map $o:G\rightarrow G\cdot f$ is separable. Example \ref{SD-type-example} in the next section gives a detailed calculation to show the orbit map is separable.\\
    (iii) For the transcendence degree, we have $\text{trdeg}_{K(G\cdot f)}K(G)=\dim G-\dim G\cdot f=\dim G(f),$ where $G(f)$ is the stabilizer of $f$ in $G$.
\end{remark}

\section{The sudden jumps and the modality}
\begin{definition}\label{xx-tpye}
    For an isolated hypersurface singularity $f \in K[[\mathbf{x}]]$, we choose the $C$-polytope $P$ as its Newton diagram (or expanding from the Newton diagram). We write the initial form as $f_0=in_P(f),\ d=v_P(f)$. We also assume that $\tau^e(f_0)<\infty$. A basis of the extended Tjurina algebra $T^e_{f_0}$ is given by $B^e_{f_0}$. 
    
    We say $f$ is of uni-deformation type if there exists a basis $B^e_{f_0}$ that satisfies the following two conditions:

    (a) for every $\mathbf{x}^{\alpha} \in B^e_{f_0}$ with $v_P(\mathbf{x}^{\alpha}) \geq d$, $v_P(\mathbf{x}^{\alpha})$ are distinct.
    
    (b)there exists $\Lambda=\{\alpha_1,\dots,\alpha_m\}$ such that $$\{\mathbf{x}^{\alpha} \in B^e_{f_0} | v_P(\mathbf{x}^{\alpha}) \geq d\}=\{\mathbf{x}^{\alpha} | \alpha \in \Lambda\}$$ and $$d \leq v_P({\mathbf{x}^{\alpha_1}})< \dots < v_P(\mathbf{x}^{\alpha_m}), \mathrm{ord}(\mathbf{x}^{\alpha_1})< \dots < \mathrm{ord}(\mathbf{x}^{\alpha_m}).$$
    
    Then we call $\Lambda$ the uni-deformation system of $f$.
\end{definition}

We recall the semicontinuity of the extended Tjurina number from \cite{Semicontinuity}.
\begin{proposition}\label{semicont-tau}
    Let $$F(\mathbf{x},\mathbf{t})=F(\mathbf{x},t_1,\dots,t_k)=f+\sum_{i=1}^k t_ig_i(\mathbf{x}) \in K[\mathbf{t}][[\mathbf{x}]],$$ where $g_i(\mathbf{x})\in K[[\mathbf{x}]]$. Then the set 
    $$U_i=\{\mathbf{t} \in K^k| \tau^e(F(\mathbf{x},\mathbf{t})) \leq i\} $$ is a Zariski open subset of $K^k$ for all $i \geq 0$.
\end{proposition}

By Proposition \ref{semicont-tau}, if we denote $$\underline{\tau}^e(f_0;g_1,\dots,g_k)=\min_{\mathbf{t}}\{\tau^e(f_0+\sum_{i=1}^k t_ig_i(\mathbf{x}))\},$$ then there exists a Zariski open subset $U \in K^k$, such that $$\tau^e(f_0+\sum_{i=1}^k t_ig_i(\mathbf{x}))=\underline{\tau}^e(f_0;g_1,\dots,g_k)$$ for any $\mathbf{t} \in U$.

Next, if we refer to general $\mathbf{t}$, it means the $\mathbf{t}$ which minimizes $\tau^e(f_0+\mathbf{t\cdot g})$.

\begin{definition}\label{sudden-jump}
    For an isolated hypersurface singularity $f \in K[[\mathbf{x}]]$ that is of uni-deformation type with the initial form $f_0=in_P(f)$ and uni-deformation system $\Lambda=\{\alpha_1,\dots,\alpha_m\}$. We denote $\underline{\tau}^e_j(f)=\underline{\tau}^e(f_0;\mathbf{x}^{\alpha_j})$ for $j \in \{1,\dots,m\}$. We call $i_1, \dots, i_k \in \{1,\dots,m-1\}$ (or $\alpha_{i_1},\dots,\alpha_{i_k}$) sudden jumps of the extended Tjurina number if 
    \begin{equation}\label{condition-sudden-jump}
        \begin{aligned}
            \underline{\tau}^e_1(f)< \dots <\underline{\tau}^e_{i_1-1}(f)<\underline{\tau}^e_{i_1+1}(f)< \dots <\underline{\tau}^e_{i_k-1}(f)<\underline{\tau}^e_{i_k+1}(f)<\underline{\tau}^e_{m}(f)\\
            <\underline{\tau}^e(f_0;\mathbf{x}^{\alpha_{i_1}},\dots,\mathbf{x}^{\alpha_{i_k}}) \leq \min \{\underline{\tau}^e_{i_1}(f), \dots, \underline{\tau}^e_{i_k}(f)\}
        \end{aligned}
    \end{equation}
    \begin{equation}\label{condition-sudden-jump-2}
        \begin{aligned}
            \mathrm{(or}\ \underline{\tau}^e_{i_1+1}(f)< \dots <\underline{\tau}^e_{i_k-1}(f)<\underline{\tau}^e_{i_k+1}(f)<\underline{\tau}^e(f_0;\mathbf{x}^{\alpha_{i_1}},\dots,\mathbf{x}^{\alpha_{i_k}}) \\ \leq \min \{\underline{\tau}^e_{i_1}(f), \dots, \underline{\tau}^e_{i_k}(f)\}\ \mathrm{if}\ i_1=1 \mathrm{)}.
        \end{aligned}
    \end{equation}
\end{definition}

\begin{example}
    As shown in figure \ref{fig-The sudden jumps}, $j=3,6$ are potential sudden jumps.
        \begin{figure}[htbp]
    \centering
    \begin{tikzpicture}[scale=2]

\coordinate (O) at (0.5,0.5);
\coordinate (A) at (1,0.75);
\coordinate (B) at (1.5,2.5);
\coordinate (C) at (2,1.25);
\coordinate (D) at (2.5,1.5);
\coordinate (E) at (3,2.5);
\coordinate (F) at (3.5,2);

\node[anchor=north] at (0.5,0) {$1$};
\node[anchor=north] at (1,0) {$2$};
\node[anchor=north] at (1.5,0) {$3$};
\node[anchor=north] at (2,0) {$4$};
\node[anchor=north] at (2.5,0) {$5$};
\node[anchor=north] at (3,0) {$6$};
\node[anchor=north] at (3.5,0) {$7$};

\draw[thick,->] (0,0) -- (4,0) node[anchor=north east]{$j$};
\draw[thick,->] (0,0,0) -- (0,3) node[anchor=north east]{$\underline{\tau}^e_j(f)$};

\draw[thick] (O) -- (A);
\draw[thick] (A) -- (B);
\draw[thick] (B) -- (C);
\draw[thick] (C) -- (D);
\draw[thick] (D) -- (E);
\draw[thick] (E) -- (F);

\foreach \point in {O, A, B, C, D, E, F}
    \fill[black] (\point) circle (1.5pt);
\end{tikzpicture}
\caption{The sudden jumps.}
\label{fig-The sudden jumps}
\end{figure}
\end{example}

\begin{example}\label{example-x3y41}
    Let $f=x^3+y^{41} \in K[[x,y]]$ with $\mathrm{char}K=5$. Then the basis of $T^e_f$ is given by
    $$B^e_f=\{1,y,y^2,\ldots,y^{40},x,xy,\ldots,xy^{39},x^2\}.$$
    Therefore, $f$ is of uni-deformation type and the uni-deformation system of $f$ is $$\Lambda=\{(1,28),(1,29),\dots,(1,39)\}.$$ $\alpha_i \in \Lambda$ is given by $\alpha_i=(1,i+27),\ i=1,2,\dots,12$.

    Through calculation, we know that $\underline{\tau}^e_{i}(f)=i+69$ for $i \neq 2,7,12$ and $\underline{\tau}^e_{2}(f)=\underline{\tau}^e_{7}(f)=\underline{\tau}^e_{12}(f)=82$. Moreover, $x^3+y^{41}+t_1xy^{29}+t_2xy^{34} \sim x^3+y^{41}+t_1xy^{29}$ for general $t_1,t_2$. In fact, apply automorphism $\phi(x)=u^{41}x,\phi(y)=u^3y$, where $u$ is a unit in $K[[x,y]]$, we have
    \begin{equation}
        \begin{aligned}
            x^3+y^{41}+t_1xy^{29}+t_2xy^{34} &\sim u^{123}(x^3+y^{41}+t_1u^5xy^{29}(1+\frac{t_2}{t_1}u^{15}y^5))\\
            &\sim x^3+y^{41}+t_1xy^{29}(u(1+(\frac{t_2}{t_1})^{\frac{1}{5}}u^3y))^5.
        \end{aligned}
    \end{equation}
    Then if suffices to find $u$ such that $u(1+(\frac{t_2}{t_1})^{\frac{1}{5}}u^3y)=1$.

    Let $F(x,y,z)=(1+z)(1+(\frac{t_2}{t_1})^{\frac{1}{5}}(1+z)^3y)-1$. We have $F(x,y,0) \in \langle x,y \rangle$ and $\frac{\partial F}{\partial z}(x,y,0) \notin \langle x,y \rangle$. By Theorem \ref{implicit-function}, there exists $z(x,y) \in \langle x,y \rangle$ such that $F(x,y,z(x,y))=0$. Choose $u=1+z(x,y)$, then $u$ is a unit and $u(1+(\frac{t_2}{t_1})^{\frac{1}{5}}u^3y)=1$.

    Therefore, $x^3+y^{41}+t_1xy^{29}+t_2xy^{34} \sim x^3+y^{41}+t_1xy^{29}$ and $\underline{\tau}^e_i(f)<\underline{\tau}^e(f_0;xy^{29},xy^{34})=\underline{\tau}^e_2(f)=\underline{\tau}^e_7(f)=82$ for $i \neq 2,7,12$. Thus, $i=2,7$ are two sudden jumps of $f$ by Definition \ref{sudden-jump}.
\end{example}

\begin{remark}
According to every example we have checked, the equation $$\underline{\tau}^e(f_0;\mathbf{x}^{\alpha_{i_1}},\dots,\mathbf{x}^{\alpha_{i_k}})=\min \{\underline{\tau}^e_{i_1}(f), \dots, \underline{\tau}^e_{i_k}(f)\}$$ always holds. Therefore, we raise a conjecture here:

\begin{conjecture}\label{strong-semicont}
    For an isolated hypersurface singularity $f \in K[[\mathbf{x}]]$, $g_1,\dots,g_k \in \mathfrak{m}^2,$ we have $\underline{\tau}^e(f;g_1,\dots,g_k)=\min_{1 \leq j \leq k} \{\underline{\tau}^e(f;g_j)\}$.
\end{conjecture}

Conjecture \ref{strong-semicont} is a more precise description of the upper semicontinuity. And if Conjecture \ref{strong-semicont} holds, we can simplify condition (\ref{condition-sudden-jump}) (or (\ref{condition-sudden-jump-2})) to
$$\underline{\tau}^e_1(f)< \dots <\underline{\tau}^e_{i_1-1}(f)<\underline{\tau}^e_{i_1+1}(f)< \dots <\underline{\tau}^e_{i_k-1}(f)<\underline{\tau}^e_{i_k+1}(f)<\underline{\tau}^e_{m}(f)
            < \min \{\underline{\tau}^e_{i_1}(f), \dots, \underline{\tau}^e_{i_k}(f)\} $$ $$\mathrm{(or}\ \underline{\tau}^e_{i_1+1}(f)< \dots <\underline{\tau}^e_{i_k-1}(f)<\underline{\tau}^e_{i_k+1}(f)< \min \{\underline{\tau}^e_{i_1}(f), \dots, \underline{\tau}^e_{i_k}(f)\}\ \mathrm{if}\ i_1=1 \mathrm{)}, $$ which is easier to check.
\end{remark}

\begin{definition}
    Let $f \in K[[\mathbf{x}]]$ be an isolated hypersurface singularity of uni-deformation type, where $P$ is the $C$-polytope given by (or expanding from) the Newton diagram of $f$. The initial form $f_0=in_P(f)$ and the uni-deformation system of $f$ is $\Lambda=\{\alpha_1,\dots,\alpha_m\}$. Assume that $i_1,\dots,i_k$ are sudden jumps of $f$. Write $l=v_P(\mathbf{x}^{\alpha_{m}})$. Choose $X=F_0/F_{l+1}$, $G=\mathcal{K}_{l}$ with the same action as defined in Section \ref{section-new-mod}.

    (a) We say $f$ is of separable deformation type (or SD type) if for every $\alpha\in \Lambda$, $g=f_0+t\mathbf{x}^{\alpha}$ with general $t$ or $t=0$, the orbit map $o:G \rightarrow G \cdot g$ is separable.

    (b) We say $f$ has separable sudden jumps if for $g=f_0+\sum_{j=1}^kt_j\mathbf{x}^{\alpha_{i_j}}$ with general $t_j$, the orbit map $o:G \rightarrow G \cdot g$ is separable.
\end{definition}

\begin{example}\label{SD-type-example}
    Consider $f=x^3+y^s \in K[[x,y]]$ with $s \geq 4,\ 3 \nmid s$, $p=\mathrm{char}K>3$. The basis of $T^e_f$ is given by $$B^e_f=\{1,y,y^2,\ldots,y^{s-1},x,xy,\ldots,xy^{s-2},x^2\}$$
    (and additionally $xy^{s-1}$ if $p \mid s$). Then $f$ is of uni-deformation type and the uni-deformation system of $f$ is $$\Lambda=\{(1,\left\lfloor {\frac{2}{3}s} \right\rfloor+1),\dots,(1,s-2)\}$$
    (and additionally $(1,s-1)$ if $p \mid s$). We will also see that $f$ is of SD type and has separable sudden jumps.

\begin{proposition}\label{example-x3ys-SD-type}
   The $f$ in Example \ref{SD-type-example} is of SD type.
\end{proposition}
\begin{proof}
        The weight vector corresponding to the Newton diagram of $f$ is $(s,3)$. We have $d=v_P(in_P(f))=3s,\ l=v_P(xy^{s-2})=4s-6$ (or $l=v_P(xy^{s-1})=4s-3$ if $p \mid s$). Therefore $X=F_{0}/F_{4s-5},\ G=\mathcal{K}_{4s-6}$ (or $X=F_{0}/F_{4s-2}$ $G=\mathcal{K}_{4s-3}$ if $p \mid s$). To simplify the discussion, we assume that $p \nmid s$. The discussion for the other case where $p \mid s$ is similar.

        Now we write $g=x^3+y^s+txy^k$ with $s \geq 4, \ 3\nmid s$, $\frac{2s}{3}<k \leq s-2$ ($k \leq s-1$ if $p \mid s$), $t$ is general or $t=0$. For $\varphi=(U,\phi) \in G$, we write 
     \begin{equation}
        \begin{aligned}
            U &=e_0+e_{10}x+e_{01}y+e_{20}x^2+e_{11}xy+e_{02}y^2+\dots,\\
            \phi(x)&=a_{10}x+a_{01}y+a_{20}x^2+a_{11}xy+a_{02}y^2+\dots,\\
            \phi(y)&=b_{10}x+b_{01}y+b_{20}x^2+b_{11}xy+b_{02}y^2+\dots.
        \end{aligned}
    \end{equation}
    Then we can write the action on $g$ in X as follows (we ignore the terms with a valuation greater than $4s-6$ or less than $3s$, and we also rewrite the symbols $e_{10},e_{01},\dots$ as $e_1,e_2,\dots$):
    \begin{equation}\label{x3ys-separable-eq}
        \begin{aligned}
            \varphi(g)&=U\cdot \phi(g)\\
            &=(e_0+e_{10}x+e_{01}y+\dots)\cdot \\
            &\bigg( (a_{10}x+a_{01}y+\dots)^3+(b_{01}y)^s+t(a_{10}x+a_{01}y+\dots)(b_{01}y)^k  \bigg).
        \end{aligned}
    \end{equation}

In fact, the term $b_{10}x$ vanishes in equation \ref{x3ys-separable-eq} since any term whose coefficient contains $b_{10}$ has a valuation greater than $v_P(b_{10}^kx^k)>l=4s-6$, which lies in $\ker \pi$, where $\pi: R \rightarrow X=F_{0}/F_{4s-5}$ is the projection map. The terms $b_{02}y^2,b_{11}xy,\dots$ vanish for the same reason.

To simplify the process, we assume that $s \equiv 1 \mod 3$ (the other case where $s \equiv 2 \mod 3$ remains the same). Continuing from \eqref{x3ys-separable-eq}, we get

\begin{equation}\label{x3ys-separable-eq-2}
        \begin{aligned}
            \varphi(g)&=U\cdot \phi(g)\\
            &=e_0a_{10}^3x^3+e_0a_{01}^3y^3+(e_{01}a_{01}^3+3e_0a_{01}^2a_{02})y^4+\dots\\
            &+\bigg(e_0\cdot (b_{01}^s+3a_{01}^2a_{0,s-2}+6a_{01}a_{02}a_{0,s-3}+\dots+ta_{0,s-k}b_{01}^k)\\
            &+e_{01}\cdot(a_{0,\frac{s-1}{3}}^3+3a_{01}^2a_{0,s-3}+\dots+ta_{0,s-k-1}b_{01}^k)+\dots\bigg)y^s\\
            &+\dots \\
            &+\bigg(e_0\cdot(6a_{10}a_{01}a_{0,\frac{2s-2}{3}}+\dots+6a_{10}a_{0,\frac{s-1}{3}}a_{0,\frac{s+2}{3}}+\dots)+\dots \bigg)xy^{\frac{2s+1}{3}}\\
            &+\dots.
        \end{aligned}
    \end{equation}
We write the coefficient of $x^iy^j$ as $c_{ij}$ in (\ref{x3ys-separable-eq-2}) to make it clear. Then the orbit map turns out to be
    
\begin{tikzcd}
&G                                    &   \longrightarrow  & G\cdot f          \\                   &{(e_0,\dots,a_{10},\dots,b_{01})}   & \mapsto & {(c_{30},c_{03},\dots,c_{0s},\dots,c_{1,s-2}).}
\end{tikzcd}\\
    Therefore the induced field extension is
    $$\widetilde{K}:=K(c_{30},c_{03},\dots,c_{0s},\dots,c_{1,s-2}) \hookrightarrow K(e_0,\dots,a_{10},\dots,b_{01}).$$

    Note that every $c_{ij}$ is a polynomial of $e_{i'j'},a_{i'j'},b_{01}$, and the degrees of $e_{i'j'},a_{i'j'}$ are less than $4$. Therefore, the degrees of minimal polynomials of $e_i,a_{ij}$ in $\widetilde{K}$ are less than $4$, thus $a_i,e_i$ are always separable elements. 

    For $b_{01}$, by calculating the dimension, we find that the stabilizer $G(g)$ of $g$ has dimension at least $1$. Therefore, the transcendence degree $trdeg_{\widetilde{K}}K(G) \geq 1$ by Remark \ref{remark-separable}.(iii). Hence we can choose $b_{01}$ as a transcendence basis, so that $K(G)$ is separably generated over $\widetilde{K}(b_{01})$, which shows that $o:G \rightarrow G \cdot g$ is separable and $f$ is of SD type as we want.
    \end{proof}

\begin{proposition}\label{example-x3ys-sudden-jump}
   The $f$ in Example \ref{SD-type-example} has separable sudden jumps.
\end{proposition}
The proof of Proposition \ref{example-x3ys-sudden-jump} will be given after Lemma \ref{x3+ys+xyk-Te-basis-lemma}.

\begin{lemma}\label{lm5.12}
    Let $I$ be an $\mathfrak{m}$-primary ideal in $K[\mathbf{x}] = K[x_1,\dots,x_n]$. Then there exists an isomorphism of \(K\)-algebras
    $$
    K[\mathbf{x}]/I \to K[[\mathbf{x}]]/IK[[\mathbf{x}]], \quad f+I \to f+ IK[[\mathbf{x}]].
    $$
    It follows that
    $$
    \dim_K (K[\mathbf{x}]/I)=\dim_K (K[[\mathbf{x}]]/IK[[\mathbf{x}]])<+\infty, 
    $$
    and as \(K\)-vector spaces, a set of representatives of a monomial basis of \(K[\mathbf{x}]/I\) necessarily gives rise to a set of representatives of a monomial basis of \(K[[\mathbf{x}]]/IK[[\mathbf{x}]]\).
\end{lemma}
\begin{proof}
Consider the natural ring homomorphism \(\varphi: K[\mathbf{x}] \to K[[\mathbf{x}]]/IK[[\mathbf{x}]]\), \(f \mapsto f+IK[[\mathbf{x}]]\). Since \(I \subseteq \ker \varphi\), it induces a homomorphism
\[
\bar{\varphi}: K[\mathbf{x}]/I \to K[[\mathbf{x}]]/IK[[\mathbf{x}]], \quad f+I \mapsto f+IK[[\mathbf{x}]].
\]

\textbf{Surjectivity:} Take any \(g \in K[[\mathbf{x}]]\) and let \(f_k\) be the truncation of \(g\) at total degree \(k-1\). Because \(I\) is \(\mathfrak{m}\)-primary, there exists \(r\) such that \(\mathfrak{m}^r \subseteq I\). For any \(k \ge r\), we have \(g-f_k \in \mathfrak{m}^k \subseteq \mathfrak{m}^r \subseteq I\), hence \(g \equiv f_k \pmod{IK[[\mathbf{x}]]}\), i.e., \(\bar{\varphi}(f_k+I)=g+IK[[\mathbf{x}]]\).

\textbf{Injectivity:} Suppose \(f \in K[\mathbf{x}]\) and \(f \in IK[[\mathbf{x}]]\), then \(f = \sum a_i h_i\) with \(a_i \in K[[\mathbf{x}]]\), \(h_i \in I\). Choose \(r\) such that \(\mathfrak{m}^r \subseteq I\) and consider the natural projection \(\pi: K[[\mathbf{x}]] \to K[[\mathbf{x}]]/(\mathfrak{m}K[[\mathbf{x}]])^r \cong K[\mathbf{x}]/\mathfrak{m}^r\). In \(K[\mathbf{x}]/\mathfrak{m}^r\) we have
\[
f+\mathfrak{m}^r = \sum (a_i+\mathfrak{m}^r)(h_i+\mathfrak{m}^r) \in I/\mathfrak{m}^r,
\]
hence \(f \in I+\mathfrak{m}^r = I\), i.e. \(\ker \bar{\varphi}=0\).

Therefore, \(\bar{\varphi}\) is a ring isomorphism and clearly \(K\)-linear, hence an isomorphism of \(K\)-algebras.
\end{proof}

\begin{proposition}\label{Calculate-basis}
    Using the notation of  Lemma \ref{lm5.12}, let $G$ be a standard basis of $I$ with respect to a given global monomial ordering $>$ (see \cite[Definitions 1.2.1, 1.2.4, 1.6.1]{singular-introduction}). For any subset $E\subset K[\mathbf{x}]$, denote by $L(E) = \langle \operatorname{LM}(\phi) \mid \phi \in E \rangle$ the leading ideal of $E$ in $K[\mathbf{x}]$, where $\operatorname{LM}(\phi)$ is the leading monomial of $\phi$ with respect to $>$ (see \cite[Definition 1.2.2]{singular-introduction}). Then the monomials in $K[\mathbf{x}] \setminus L(G)$ constitute a $K$-basis for the vector space $K[[\mathbf{x}]]/IK[[\mathbf{x}]]$.
\end{proposition}

\begin{proof}
    Combining  Lemma \ref{lm5.12} and \cite[Corollary 7.5.6]{singular-introduction} immediately yields that the monomials in $K[\mathbf{x}] \setminus L(I)$ constitute a $K$-basis for the vector space $K[[\mathbf{x}]]/IK[[\mathbf{x}]]$. By the definition of a standard basis (see \cite[Definition 1.6.1]{singular-introduction}), we have $L(G) = L(I)$, completing the proof.
\end{proof}

\begin{lemma}\label{x3+ys+xyk-Te-basis-lemma}
    For $t\neq 0$, let $g=x^3+y^s+txy^k$ with $s \geq 4$, $2s<3k<3s$ as in Proposition \ref{example-x3ys-SD-type}. Then $T^e_g$ has a monomial basis given by 
    $$
    \begin{cases}
    \{1,y,y^2,\ldots,y^{s-1},x,xy,\ldots,xy^{k-1},x^2\}, & p\nmid 3k-2s \\
    \{1,y,y^2,\ldots,y^{s-1},x,xy,\ldots,xy^{s-2},x^2\}, & p\mid 3k-2s,p\nmid k,p\nmid s \\
    \{1,y,y^2,\ldots,y^{s-1},x,xy,\ldots,xy^{s-1},x^2\}, & p\mid k,p\mid s
    \end{cases}
    $$
    or equivalently,
    $$
    \begin{cases}
    \{1,y,y^2,\ldots,y^{s-1},x,xy,\ldots,xy^{k-1},x^2\}, & p\nmid 3k-2s \\
    \{1,y,y^2,\ldots,y^{2s-k-2},x,xy,\ldots,xy^{k-1},x^2\}, & p\mid 3k-2s,p\nmid k,p\nmid s \\
    \{1,y,y^2,\ldots,y^{2s-k-1},x,xy,\ldots,xy^{k-1},x^2\}, & p\mid k,p\mid s
    \end{cases},
    $$
    where $p=\operatorname{char}K$ is either zero or a prime number greater than $3$. Therefore, the extended Tjurina number is
    $$\tau^e(g)=\dim_K T^e_g  =\begin{cases}
    k+s+1, & p\nmid 3k-2s \\
    2s, & p\mid 3k-2s,p\nmid k,p\nmid s \\
    2s+1, & p\mid k,p\mid s.
    \end{cases} $$
\end{lemma}

\begin{proof}
Since $s\ge k+1$ and $3k\ge 2s+1$, we have $2k-s\ge  2$. These three inequalities will be tacitly employed throughout the subsequent proof without explicit mention.

\textbf{(I)} $p\nmid 3k-2s$.

Since
$$
\begin{pmatrix} g \\ xg_x \\ yg_y \end{pmatrix}
= \begin{pmatrix} 1 & 1 & t \\ 3 & 0 & t \\ 0 & s & kt \end{pmatrix}
\begin{pmatrix} x^3 \\ xy^{k} \\ y^s \end{pmatrix},
$$
and the determinant of the coefficient matrix is \(t(2s-3k) \neq 0\) for \(t \neq 0\), it follows that
$$
x^3,\; xy^{k},\; y^s \in \langle g,\; xg_x,\; yg_y\rangle.
$$
Consequently,
$$
x^3,\; 3x^2y+ty^{k+1},\; xy^{k},\; y^s \in \langle g,\; \mathfrak{m}\cdot j(g) \rangle.
$$
Furthermore,
\[
\begin{pmatrix} 
  g \\ 
  xg_x \\ 
  yg_x \\ 
  xg_y \\ 
  yg_y 
\end{pmatrix}
= 
\begin{pmatrix} 
    1 & 0 & t & 1  \\ 
    3 & 0 & t & 0  \\ 
    0 & 1 & 0 & 0  \\  
    0 & \frac{kt}{3}y^{k-2} & sy^{s-k-1} & -\frac{kt^2}{3}y^{2k-s-1}  \\  
    0 & 0 & kt & s  \\
\end{pmatrix}
\begin{pmatrix} 
  x^3 \\ 
  3x^2y+ty^{k+1} \\ 
  xy^{k} \\ 
  y^s 
\end{pmatrix}.
\]
Therefore,
$$
\langle g,\mathfrak{m}\cdot j(g) \rangle= \langle x^3,3x^2y+ty^{k+1},xy^k,y^s \rangle.
$$
Using the lexicographical ordering $>_{lp}$ (see \cite[Example 1.2.8]{singular-introduction}), one can apply Buchberger's criterion (see \cite[Theorem 1.7.3]{singular-introduction}) to verify that $G=\{ x^3,\ 3x^2y + ty^{k+1},\ xy^k,\ y^s \}$ is a standard basis and $L(G)= \langle x^3,x^2y,xy^k,y^s \rangle$.
Therefore, by Proposition \ref{Calculate-basis}, 
$T_g^e=K[[x,y]]/\langle g,\mathfrak{m}\cdot j(g) \rangle$ has a $K$-basis represented by 
$$
\{1,y,y^2,\ldots,y^{s-1},x,xy,\ldots,xy^{k-1},x^2\}.
$$

\textbf{(II)} $p\mid 3k-2s$ with $p\nmid k,\ p\nmid s$.

Since $g=\frac{1}{3}xg_x+\frac{2}{3k}yg_y$ , we have $\langle g, \mathfrak{m}\cdot j(g) \rangle= \langle xg_x,yg_x,xg_y,yg_y \rangle$. Note that
$$
\begin{aligned}
\left(1+\frac{k^2t^3}{3s^2}y^{3k-2s}\right)y^{2s-k-1}
=&\frac{k^2t^2}{3s^2}y^{k-1}(3x^2+ty^k)+\left(-\frac{kt}{s^2}x+\frac{1}{s}y^{s-k}\right)(ktxy^{k-1}+sy^{s-1})\\
=&\frac{k^2t^2}{3s^2}y^{k-1}g_x+\left(-\frac{kt}{s^2}x+\frac{1}{s}y^{s-k}\right)g_y.
\end{aligned}
$$
It follows that
\[
\begin{pmatrix} 
  3x^3+txy^k \\ 
  3x^2y+ty^{k+1} \\ 
  2txy^{k}+3y^s \\ 
  y^{2s-k-1} 
\end{pmatrix}
=
\begin{pmatrix} 
  1 & 0 & 0 & 0 \\ 
  0 & 1 & 0 & 0 \\ 
  0 & 0 & 0 & \frac{2}{k} \\ 
  0 & \frac{\frac{k^2t^2}{3s^2}y^{k-2}}{1+\frac{k^2t^3}{3s^2}y^{3k-2s}} & \frac{\frac{-kt}{s^2}}{1+\frac{k^2t^3}{3s^2}y^{3k-2s}} & \frac{\frac{1}{s}y^{s-k-1}}{1+\frac{k^2t^3}{3s^2}y^{3k-2s}}
\end{pmatrix}
\begin{pmatrix} 
  xg_x \\ 
  yg_x \\ 
  xg_y \\
  yg_y
\end{pmatrix}
\]
where the determinant of the coefficient matrix is
$$\frac{\frac{2t}{s^2}}{1+\frac{k^2t^3}{3s^2}y^{3k-2s}}
=\frac{\frac{2t}{s^2}}{1+\frac{4t^3}{27}y^{3k-2s}}
\neq 0\quad \text{for $t\neq 0$}.$$
Therefore,
$$
\langle g,\mathfrak{m}\cdot j(g) \rangle
= 
\langle 3x^3+txy^k,3x^2y+ty^{k+1},2txy^{k}+3y^s,y^{2s-k-1} \rangle.
$$
Proceeding as in the first case, one finds that this set of elements constitutes a standard basis and that the quotient $T_g^e= K[[x,y]]/\langle g, \mathfrak{m} \cdot j(g) \rangle$ admits a $K$-basis represented by
$$
\{1,y,y^2,\ldots,y^{s-1},x,xy,\ldots,xy^{s-2},x^2\}.
$$
Choosing other global orderings may yield different but equivalent sets of monomial basis representatives.

\textbf{(III)} $p\mid k$ and $p\mid s$. 
Then $g_y=ktxy^{k-1}+sy^{s-1}=0$ and $\langle g,\mathfrak{m}\cdot j(g) \rangle= \langle g,xg_x,yg_x \rangle$. Note that
$$
\begin{aligned}
&\left(1+\frac{4}{27}t^3y^{3k-2s}\right)y^{2s-k} \\
=&\left(-\frac{2}{3}tx+y^{s-k}\right)(x^3+y^s+txy^k)+\left(\frac{2}{9}tx^2-\frac{1}{3}xy^{s-k}+\frac{4}{27}t^2y^k\right)(3x^2+ty^k) \\
=&\left(-\frac{2}{3}tx+y^{s-k}\right)g+\left(\frac{2}{9}tx^2-\frac{1}{3}xy^{s-k}+\frac{4}{27}t^2y^k\right)g_x. \\
\end{aligned}
$$
It follows that
\[
\begin{pmatrix} 
  3x^3+txy^k \\ 
  3x^2y+ty^{k+1} \\ 
  2txy^{k}+3y^s \\ 
  y^{2s-k} 
\end{pmatrix}
=
\begin{pmatrix} 
  0 & 1 & 0 \\ 
  0 & 0 & 1  \\ 
  3 & -1 & 0 \\ 
  \frac{-\frac{2}{3}tx+y^{s-k}}{1+\frac{4}{27}t^3y^{3k-2s}} & \frac{\frac{2}{9}tx-\frac{1}{3}y^{s-k}}{1+\frac{4}{27}t^3y^{3k-2s}} & \frac{\frac{4}{27}t^2y^{k-1}}{1+\frac{4}{27}t^3y^{3k-2s}}  
\end{pmatrix}
\begin{pmatrix} 
  g \\ 
  xg_x \\ 
  yg_x \\ 
\end{pmatrix}
\]
and
\[
\begin{pmatrix} 
  g \\ 
  xg_x \\ 
  yg_x \\ 
\end{pmatrix}
= 
\begin{pmatrix} 
    \frac{1}{3} & 0 & \frac{1}{3} & 0  \\ 
    1 & 0 & 0 & 0  \\ 
    0 & 1 & 0 & 0  \\  
\end{pmatrix}
\begin{pmatrix} 
  3x^3+txy^k \\ 
  3x^2y+ty^{k+1} \\ 
  2txy^{k}+3y^s \\ 
  y^{2s-k} 
\end{pmatrix}.
\]
Therefore,
$$
\langle g,\mathfrak{m}\cdot j(g) \rangle
=
\langle 3x^3+txy^k,3x^2y+ty^{k+1},2txy^{k}+3y^s,y^{2s-k} \rangle.
$$
Proceeding as in the first case, one finds that this set of elements constitutes a standard basis and that the quotient $T_g^e= K[[x,y]]/\langle g, \mathfrak{m} \cdot j(g) \rangle$ admits a $K$-basis represented by
$$
\{1,y,y^2,\ldots,y^{2s-k-2},x,xy,\ldots,xy^{k-1},x^2\}.
$$
Choosing other global orderings may yield different but equivalent sets of monomial basis representatives.
\end{proof}

\textit{Proof(of Proposition \ref{example-x3ys-sudden-jump}).}
By Lemma \ref{x3+ys+xyk-Te-basis-lemma}, we can see that the sudden jump of $f$ is $(1,k)$ for $p \mid 3k-2s, \frac{2}{3}s \leq k \leq s-2$ (or $\frac{2}{3}s \leq k \leq s-1$ if $p \mid s$). Similar to Example \ref{example-x3y41}, we can show $$x^3+y^s+\sum_{(1,k)\mathrm{\ is\ sudden\ jump}}t_kxy^k$$ is in the same orbit with $x^3+y^s+t_{k_0}xy^{k_0}$, where $(1,k_0)$ is the minimal sudden jump. Therefore, by Proposition \ref{example-x3ys-SD-type}, the orbit map of $x^3+y^s+t_{k_0}xy^{k_0}$ is separable and $f$ has separable sudden jumps.
\qed
    
\end{example}

Then we can state the theorem. In the following of this section, $f$ always refers to an isolated hypersurface singularity of uni-deformation type and $\tau^e(in_P(f))< \infty$, where $P$ is the $C$-polytope given by (or expanding from) the Newton diagram of $f$.

\begin{theorem}\label{main-thm}
    If $f$ is of SD type and $f$ has separable sudden jumps, then the $\mathcal{K}$-modality of $f$ is greater than or equal to the number of sudden jumps of $f$.
\end{theorem}

\begin{proof}
    Write $f_0=in_P(f)$ and $d=v_P(f_0).$ Let $\Lambda=\{\alpha_1,\dots,\alpha_m\}$ be the uni-deformation system of $f$. Assume $i_1,\dots,i_k$ are sudden jumps of $f$.

    Write $l=v_P(\mathbf{x}^{\alpha_{m}})$. Choose $X=F_0/F_{l+1}$, $G=\mathcal{K}_{l}$ with the same action as defined in Section \ref{section-new-mod}. By Definition \ref{def-mod}, there exists an open neighborhood $U \subset X$ of $f$, such that $$G\text{-mod}(f)=G\text{-mod}(U)=\max_{i \geq 0}\{\dim U(\leq i)-i\}.$$

    We denote $g_0=f_0+t_m\mathbf{x}^{\alpha_m}$ for general $t_m \in K$ and denote $i_0=\dim U \cap G \cdot g_0$. Since $f$ has SD-type, the orbit map of $g_0$ is separable. We have
    $$\dim G\cdot g_0=\dim \widetilde{T}_{g_0}(Gg_0)=\dim X- \underline{\tau}^e(g_0)$$ by Corollary \ref{dim-Gf-cor}.

    On the other hand, let $U_1=\{f_0+\sum_{j=1}^kt_j\mathbf{x}^{\alpha_{i_j}}\mid t_j \in K^{\times} \mathrm{\ are\ general}\}$. Since $f$ has separable sudden jumps, each orbit map of $g\in U_1$ is separable, and therefore $$\dim G\cdot g=\dim \widetilde{T}_{g}(Gg)=\dim X- \underline{\tau}^e(g).$$
    Since $\underline{\tau}^e(g_0)<\underline{\tau}^e(g)$, we have $\dim G\cdot g_0>\dim G \cdot g$. Therefore, $U \cap U_1 \subset U(\leq i_0)$ and $U_1 \cap G \cdot g_0=\emptyset$.

    We claim that: $\mathbf{x}^{\alpha_{i_j}} \notin {T}_{g_0}(G \cdot g_0)$ for $j=1,\dots,k$.

    Since $o:G \rightarrow G \cdot g_0$ is separable, ${T}_{g_0}(G \cdot g_0)=\langle g_0 \rangle+\mathfrak{m}\cdot \langle \frac{\partial g_0}{\partial x_1},\dots,\frac{\partial g_0}{\partial x_n}\rangle+F_{l+1}/F_{l+1}$. If $\mathbf{x}^{\alpha_{i_j}} \in {T}_{g_0}(G \cdot g_0)$, we have
    \begin{equation}\label{main-contradiction}
        \mathbf{x}^{\alpha_{i_j}}=u_0f_0+u_0\mathbf{x}^{\alpha_m}+\sum_{t=1}^n u_t(\frac{\partial f_0}{\partial x_t}+\frac{\partial \mathbf{x}^{\alpha_m}}{\partial x_t}),\ u_t \in \mathfrak{m} \mathrm{\ for\ }t=1,\dots,n.
    \end{equation}
    But $\mathbf{x}^{\alpha_j} \notin T_{f_0}(G \cdot f)=\langle f_0 \rangle+\mathfrak{m}\cdot \langle \frac{\partial f_0}{\partial x_1},\dots,\frac{\partial f_0}{\partial x_n}\rangle$ and $u_t \frac{\partial x^{\alpha_m}}{\partial x_t}$ cannot generate $\mathbf{x}^{\alpha_{i_j}}$ because $\mathrm{ord}(\mathfrak{m}\cdot \frac{\partial \mathbf{x}^{\alpha_m}}{\partial x_t}) \geq \mathrm{ord}(\mathbf{x}^{\alpha_m})>\mathrm{ord}(\mathbf{x}^{\alpha_{i_j}}).$ Then equation (\ref{main-contradiction}) cannot hold and we prove the claim. 

    Therefore we have $$T_{g_0}(U \cap G \cdot g_0) \oplus span \langle \mathbf{x}^{\alpha_{i_1}},\dots,\mathbf{x}^{\alpha_{i_k}} \rangle \subset T_{g_0}U(\leq i_0)$$ and thus $$G-mod(f) \geq \dim U(\leq i_0)-i_0=\dim U(\leq i_0)-\dim U \cap G\cdot g_0 \geq k.$$
\end{proof}

\section{The classification in characteristic $p >3$}\label{classification p}
We now state our classification results.
\begin{proposition}\label{class-xy}
    The following hypersurface singularities are the only candidates for modality 1 in $K[[x,y]]$:
\renewcommand\arraystretch{1.5}

\begin{longtable}{|c|c|c|}
\caption{}\label{table-xy}
\\
\hline
Symbol&Form& condition  \\
\hline
$E_{0,s}$&$x^3+y^{s}$&$s \geq 6$ \\
\hline
$E_{r,0}$&$x^3+xy^{r}$&$r\geq 4$ \\
\hline
$E_{k,s}$&$x^3+y^s+xy^k$& $k\geq 3,s\geq 4,p\mid 3k-2s$\\
\hline
$E_{k,s,l}$&$x^3+y^s+\lambda xy^k+xy^l$& $k\geq 3,s\geq 4,l>k,p\mid 3k-2s,p \nmid 3l-2s,\lambda \in K$\\
\hline
$E_{2t,3t,0}$&$x^3+xy^{2t}+\lambda y^{3t}$& $t\geq 2,\lambda \neq 0,p \neq 31$\\
\hline
$E_{2t,3t,l}$&$x^3+xy^{2t}+\lambda y^{3t}+y^l$& $t\geq 2,l>3t,\lambda \neq 0$\\
\hline
$W_{12}$&$x^4+y^5$ &$p \neq 5$ \\
\hline
$W_{12}'$&$x^4+y^5+x^2y^3$&$p \neq 5$ \\
\hline
$W_{13}$&$x^4+xy^4$&\\
\hline
$W_{13}'$&$x^4+xy^4+y^6$ &\\
\hline
$W_{1,0}$&$x^4+x^2y^3+\lambda y^6$& $\lambda \neq 0,\frac{1}{4}$  \\
\hline
$W_{1,0}'$&$x^4+x^2y^3+\lambda y^6+y^7$& $\lambda \neq 0,\frac{1}{4}$ \\
\hline
$W_{1,t}$&$x^4+x^2y^3+y^t$& $t \geq 7$\\
\hline
$W_{1,0}^\#$&$x^4+y^6$&\\
\hline
$W_{1,0}^{\#'}$&$x^4+x^2y^4+y^6$&\\
\hline
$W_{17}$&$x^4+xy^5$&$p \neq 5$\\
\hline
$W_{17}'$&$x^4+xy^5+y^7$&$p \neq 5$\\
\hline
$W_{17}''$&$x^4+xy^5+y^8$&$p \neq 5$\\
\hline
$W_{18}$&$x^4+y^7$&$p \neq 7$\\
\hline
$W_{18}'$&$x^4+y^7+x^2y^4$&$p \neq 7$\\
\hline
$W_{18}$&$x^4+y^7+x^2y^5$&$p \neq 7$\\
\hline
$Z_{6m+5}$&$x^3y+y^{3m+2}$&$m \geq 1$\\
\hline
$Z_{6m+6}$&$x^3y+xy^{2m+2}$&$m \geq 1$\\
\hline
$Z_{6m+7}$&$x^3y+y^{3m+3}$&$m \geq 1$\\
\hline
$Z_{k,s,l}$&$x^3y+y^s+\lambda xy^k+xy^l$&$k\geq 4,s \geq 5,l>k,p\mid 3k-2s-1,p \nmid 3l-2s-1,\lambda \in K$\\
\hline
$Z_{k,s,l}'$&$x^3y+xy^{2t+1}+\lambda y^{3t+1}+y^l$& $t\geq 2,l>3t+1,p \nmid l-3t-1,\lambda \neq 0$\\
\hline
$T_{4,s,2}$&$x^4+x^2y^2+y^s$&$s \geq 5$\\
\hline
$T_{r,s,2}$&$x^r+x^2y^2+y^s$&$r,s \geq 5$\\
\hline
$T_{4,4,2}$&$x^4+\lambda x^2y^2+y^4$&$\lambda^2 \neq 4$\\
\hline
\end{longtable}

\end{proposition}

\begin{proposition}\label{class-xyz}
    The following hypersurface singularities are the only candidates for modality 1 in $K[[x,y,z]]$:
\renewcommand\arraystretch{1.5}
\begin{longtable}{|c|c|c|}
\caption{}\label{table-xyz}
\\
\hline
Symbol&Form& condition  \\
\hline
$T_{3,3,3}$&$x^3+y^3+z^3+\lambda xyz $&$\lambda^3+27\neq 0$ \\
\hline
$T_{r,s,t}$&$x^r+y^s+z^t+xyz$&$\max \{r,s,t\}\geq 4$ \\
\hline
$Q_{6m+4}$&$x^3+yz^2+y^{3m+1}$&$m\geq1$ \\
\hline
$Q_{6m+5}$&$x^3+yz^2+xy^{2m+1}$& $m \geq 1$ \\
\hline
$Q_{6m+6}$&$x^3+yz^2+y^{3m+2}$&$m\geq1$ \\
\hline
$Q_{k,s,l}$&$x^3+yz^2+y^s+\lambda xy^k+xy^l$& $k\geq 3,s\geq 4,l>k,p\mid 3k-2s,p \nmid 3l-2s,\lambda \in K$\\
\hline
$Q_{r,s,l}'$&$x^3+yz^2+xy^{2t}+\lambda y^{3t}+y^l$& $t\geq 2,l>3t,p \nmid l-3t,\lambda \neq 0$\\
\hline
$S_{11}$&$x^2z+yz^2+y^4$ & \\
\hline
$S_{11}'$&$x^2z+yz^2+y^4+\lambda x^2y^2$ & \\
\hline
$S_{12}$&$x^2z+yz^2+xy^3$& \\
\hline
$S_{1,0}$&$x^2z+yz^2+x^2y^2+\lambda y^5$&$\lambda \neq 0$\\
\hline
$S_{1,0}^{1}$&$x^2z+yz^2+x^2y^2+\lambda y^5+y^6$&$\lambda \neq 0$\\
\hline
$S_{1,0}^{2}$&$x^2z+yz^2+x^2y^2+xy^4$&\\
\hline
$S_{1,0}^{3}$&$x^2z+yz^2+y^5$ &$p \neq 5$\\
\hline
$S_{1,0}^{4}$&$x^2z+yz^2+x^2y^3+y^5$ &$p \neq 5$\\

\hline
$S_{1,0,t}$&$x^2z+yz^2+x^2y^2+ y^t$&$6\leq t<s+2$\\
\hline
$S_{1,s,0}$&$x^2z+yz^2+x^2y^2+ xy^s$&$t\geq 2s-2$\\
\hline
$S_{1,s,t}$&$x^2z+yz^2+x^2y^2+ xy^s+\lambda y^{t}$&$s\geq 5,s+2\leq t\leq 2s-3,\lambda \neq 0$\\
\hline
$S_{16}$&$x^2z+yz^2+xy^4$&\\
\hline
$S_{16}'$&$x^2z+yz^2+xy^4+y^6$&\\
\hline
$S_{16}''$&$x^2z+yz^2+xy^4+y^7$&\\
\hline
$S_{17}$&$x^2z+yz^2+y^6$ &\\
\hline
$S_{17}'$&$x^2z+yz^2+y^6+x^2y^3$ &\\
\hline
$S_{17}''$&$x^2z+yz^2+y^6+x^2y^4$ &\\
\hline
$U_{12}$&$x^3+xz^2+y^4$&  \\
\hline
$U_{12}'$&$x^3+xz^2+y^4+x^2y^2$&  \\
\hline
$U_{1,0}$&$x^3+xz^2+xy^3+\lambda y^3z$& $\lambda^2 \neq 0,-1$  \\
\hline
$U_{1,0}'$&$x^3+xz^2+xy^3+\lambda y^3z+y^4z$& $\lambda^2 \neq 0,-1$  \\
\hline
$U_{1,t}$&$x^3+xz^2+xy^3+y^tz$& $t \geq 4$  \\
\hline
$U_{16}$&$x^3+xz^2+y^5$& $p \neq 5$  \\
\hline
$U_{16}'$&$x^3+xz^2+y^5+x^2y^3$& $p \neq 5$  \\
\hline
$U_{*}$&$x^3+xz^2+y^3z$&   \\
\hline
$U_{*}'$&$x^3+xz^2+y^3z+xy^4$&   \\
\hline
\end{longtable}
\end{proposition}

\begin{proposition}\label{class-n>3}
    All unimodal hypersurface singularities in $K[[x_1,\dots,x_n]]$ with $n \geq 4$ must be of the form $g(x_1,x_2)+x_3^2+\dots+x_n^2$ or $h(x_1,x_2,x_3)+x_4^2+\dots+x_n^2$, where $g$ (resp. $h$) is one of the forms in Table \ref{table-xy} (resp. Table \ref{table-xyz}).
\end{proposition}

We begin with $n=2$. 
\subsection{Unimodal hypersurface singularities in $K[[x,y]]$}
Assume $l=\mathrm{ord}(f) \geq 2$. By Proposition \ref{nl-bound}, we have $l \leq 4$.

If $l=2$, we have the following splitting lemma for char$K \neq 2$ from \cite{right-simple}.

Let $f\in K[[{\bf x}]]=K[[x_1,\ldots,x_n]]$. We denote by
$$H(f):=\big( \frac{\partial^2 f}{\partial x_i\partial x_j}(0)\big)_{i,j=1,\ldots,n}\in \mathrm{Mat}(n\times n, K)$$
the {\em Hessian (matrix)} of $f$ and by $\mathrm{crk}(f):=n-\mathrm{rank}(H(f))$ the {\em corank} of $f$. 

\begin{lemma} \label{splitting}
If $f\in \mathfrak{m}^2\subset K[[{\bf x}]], \mathrm{char}(K)>2,$ has corank $\mathrm{crk}(f)=k\geq 0$, then 
$$f\sim g(x_1,\ldots, x_k)+x_{k+1}^2+\ldots + x_{n}^2$$ 
with $g\in \mathfrak{m}^3$. 
\end{lemma}

Using Lemma \ref{splitting}, we can see that for $f \in \mathfrak{m}^2\subset K[[x,y]]$ with $\mathrm{ord}(f)=2$, then $f$ must be contact equivalent to $A_k:\ x^2+y^{k+1},\ k \geq 1$, which is simple.

Now assume $ord(f)=3$. Then $j_3(f)$ has one of the following forms: $x^3,x^2y,x^2y+xy^2$.

We will provide a detailed classification procedure for $f$ with $j_3(f)=x^3$. 

\begin{proposition}\label{x^3}
    If $j_3(f) \sim x^3,$ then $f$ belongs to the family $E$.
\end{proposition}
\begin{proof}
    Write $g=x^3$. Then $\widetilde{T}_g(\mathcal{K}g)=g+\mathfrak{m}\cdot j(g)=  \langle x^3,x^2y \rangle$. Therefore for any $l \geq 4$, we can find $$C=\mathrm{span}\langle xy^i,y^j \mid 3 \leq i \leq l-1, 4 \leq j \leq l \rangle$$ such that $$P_{3,l} \subset C+\widetilde{T}_g(\mathcal{K}_{l}g)\cap P_{3,l}.$$ By Theorem \ref{complete-cor}, we have
    \begin{equation}\label{x3-eq}
        f \sim x^3+\sum_i a_ixy^i+\sum_j b_jy^j=x^3+a(y)xy^r+b(y)y^s
    \end{equation}
    for some $r\geq 3,s \geq 4$, $a(y),b(y)$ are either units or $0$.

    If $a(y)=b(y)=0,$ then $f \sim x^3$, which is not isolated. If $a(y)=0$ and $b(y)$ is a unit, then $f \sim x^3+b(y)y^s$. Apply the automorphism $\phi(x)=b(y)^{\frac{1}{3}}x,\phi(y)=y$. Then $$f \sim b(y)(x^3+y^s) \sim x^3+y^s.$$ If $b(y)=0$ and $a(y)$ is a unit, similarly we have $$f \sim x^3+xy^r.$$

    Next we assume that both $a(y)$ and $b(y)$ are units. Then $f$ is convenient. The Newton diagram depends on $r,s$.

    \textbf{(I)} If $2s<3r$, then $in_f(f)=x^3+y^s$. The weight vector corresponding to the Newton diagram $P$ is $(s,3)$ and $d=v_f(f)=3s$. 
    
    If $3 \nmid s$, as we have seen in Example \ref{SD-type-example}, the uni-deformation system of $f$ is given by $\Lambda=\{(1,\left\lfloor {\frac{2}{3}s} \right\rfloor+1),\dots,(1,s-2)\}$ (and additionally $(1,s-1)$ if $p \mid s$). $f$ is of SD type and has separable sudden jump $(1,k)$ with $p \mid 2s-3k, \left\lfloor {\frac{2}{3}s} \right\rfloor+1 \leq k \leq s-3$. Thus, if there exists a $k$ such that $$\left\lfloor {\frac{2}{3}s} \right\rfloor+1 \leq k \leq k+p \leq s-3$$ and $p\mid 3k-2s$, we have $\mathcal{K}\text{-mod}(f)\geq 2$ by Theorem \ref{main-thm}.

    If $3 \mid s$, we can use the same method as Example \ref{SD-type-example} to show $f$ is of SD type. We can also find that $(1,k)$ with $p \mid 2s-3k, \frac{2}{3}s \leq k \leq s-3$ are separable sudden jumps.
    
    Therefore, the unimodal hypersurface singularity $f \in K[[x,y]]$ with $in_f(f)=x^3+y^s$ must satisfy the condition that there exists no $k$ such that 
    \begin{equation}\label{condition-k+p}
        \left\lfloor {\frac{2}{3}s} \right\rfloor+1 \leq k < k+p \leq s-3 \text{ and } p\mid 3k-2s.
    \end{equation}
    
    Since the regular bases ${\bf x}^\alpha$ of $in_f(f)$ with $v_f({\bf x}^\alpha)>3s$ are $$\{xy^{\left\lfloor {\frac{2}{3}s} \right\rfloor+1},\dots, xy^{s-2}\}\ (resp.\ \{xy^{\left\lfloor {\frac{2}{3}s} \right\rfloor+1},\dots, xy^{s-i}\} \text{ if }p\mid s ),$$ we know that $$f \sim x^3+y^s+\sum_{l \geq k}c_lxy^l,$$ where $k \geq \left\lfloor {\frac{2}{3}s} \right\rfloor+1, c_l \in K$ and $l \leq s-2$ (resp. $l \leq s-1$ if $p \mid s$) by Corollary \ref{normal-cor}. If $c_k=0$, then $f \sim x^3+y^s$. Next we assume $c_k \neq 0$. Therefore, $f \sim x^3+y^s+e(y)\cdot xy^k$, where $e(y)=c_k+c_{k+1}y+\dots$ is a unit of $R$. We rewrite $e(y)=\sum_{i \geq 0}e_iy^i $ for convenience.
    
    If $p \nmid 3k-2s$, then $f \sim x^3+y^s+xy^k$. In fact, consider the function 
    $$F(z)=z^{3k-2s}\sum_{i \geq 0}e_iy^iz^{3i}-e_0.$$
    We have $F(1)\in \langle y \rangle K[[y]]$, and $$\ F'[1]=(3k-2s)\sum_{i \geq 0}e_iy^i-3\sum_{i \geq 1}i e_iy^i$$ is a unit since $3k-2s \neq 0 \mathrm{\ and\ } p \nmid 3k-2s$. Apply Theorem \ref{implicit-function} to the function $G(z)=F(z+1)$, there exists a $\widetilde{z}(y)$ such that $G(\widetilde{z}(y))=0.$ Let $z(y)=\widetilde{z}(y)+1$, then $z(y)$ is a unit and $F(z(y))=0$, that is, $z(y)^{3k-2s}e(z(y)^3y)=e_0$.

    Using the automorphism $\phi(x)=z(y)^k x$ and $\phi(y)=z(y)^3 y$, we have $$f \sim z(y)^{3s}(x^3+y^s+z(y)^{3k-2s}e(z(y)^3y)xy^k) \sim x^3+y^s+e_0xy^k.$$ Then apply $\xi (x)=\alpha x,\ \xi(y)=\beta y$ with $\alpha,\ \beta \in F$ satisfying $\alpha^3=e_0\alpha\beta^k,\ \alpha\beta^r=\beta^s$ (such $\alpha,\ \beta$ exists since $3k-2s \neq 0$), we have $$f \sim x^3+y^s+xy^k\in E_{0,s,k}.$$ We call the method we use here the $\alpha,\beta$-trick, which is similar to what we used in Example \ref{SD-type-example}.

    If $p \mid 3k-2s$, choose $l$ to be the smallest $l$ that satisfies $c_l \neq 0, l \leq s-2$ (resp. $l \leq s-1$ if $p \mid s$) and $p \nmid 3l-2s$. If such $l$ does not exist, then $f=x^3+y^s+c_kxy^k \sim x^3+y^s+xy^k \in E_{k,s}$. If such $l$ exists, then $l$ satisfies $k<l<k+p$, otherwise $k<k+p\leq l-1\leq s-3$, which contradicts condition \eqref{condition-k+p} (if moreover $p \mid s$, then $p \mid k$, which means that $k+p<s-3$ still holds since $p>3$, leading to the same contradiction). Now we write $f=x^3+y^s+xy^k+xy^l\cdot e'(y)$. Using the same technique as the implicit function theorem (working on the terms $x^3,y^s,xy^l$), we get
    $$f \sim x^3+y^s+\widetilde{e}(y)^{3k-2s}xy^k+xy^l,$$ where $\widetilde{e}(y) \in R$ is another unit. Since $p \mid 3k-2s$, we can write $\widetilde{e}(y)^{3k-2s}$ as $$\widetilde{e}(y)^{3k-2s}=\widetilde{e}_0+\widetilde{e}_1y^p+\widetilde{e}_2y^{2p}+\dots.$$ Then $f \sim x^3+y^s+\widetilde{e}_0xy^k+xy^l+\widetilde{e_1}xy^{k+p}+\dots$.

    If $k+p \geq s-1$, then $$v_f(xy^{k+p})>\max_{\mathbf{x}^\alpha \mathrm{\ is\ a\ regular\ basis}} \{v_f(\mathbf{x}^\alpha)\}=v_f(xy^{s-2})$$ (if additionally $p \mid s$, then $p \mid k$ hence $k+p \geq s$, we still have $v_f(xy^{k+p})>v_f(xy^{s-1})$). By Corollary \ref{normal-deter}, we have $$f \sim x^3+y^s+\widetilde{e}_0xy^k+xy^l\in E_{k,s,l}.$$

    The last case to deal with is $k+p=s-2$ (otherwise $k+p \leq s-3$ will lead to a contradiction to condition \eqref{condition-k+p}). Note that in this case $p \mid s$ and $p^2 \mid 2k-2s$ cannot occur. Now we have 
    \begin{equation}
        \begin{aligned}
        f &\sim x^3+y^s+\widetilde{e}_0xy^k+xy^l+\widetilde{e_1}xy^{k+p}\\
        & \sim x^3+y^s+\widetilde{e}_0xy^k+xy^l(1+\widetilde{e}_1xy^{k+p-l}).
    \end{aligned}
    \end{equation}
    Applying the automorphism $\phi(x)=z(y)^sx,\phi(y)=z(y)^3y$, we have
    \begin{equation}\label{eq-k+p}
        f \sim z(y)^{3s}(x^3+y^s+\widetilde{e}_0z(y)^{3k-2s}xy^k+(1+\widetilde{e}_1z(y)^{k+p-l}xy^{k+p-l})z(y)^{3l-2s}xy^l.
    \end{equation}
    We hope to choose a suitable $z(y)$ such that the $xy^{k+p}$ term vanishes in \eqref{eq-k+p}.

    By Corollary \ref{normal-deter}, we can ignore all terms of the form $xy^i$ with $i>k+p$. Therefore, we can assume $z(y)=1+ty$ and apply the method of undetermined coefficients to find $t \in K$. Then \eqref{eq-k+p} becomes
    \begin{equation}
        \begin{aligned}
            f &\sim x^3+y^s+\widetilde{e}_0(1+\frac{3k-2s}{p}t^{p}y^{p})xy^k+(1+\widetilde{e}_1xy^{k+p-l})xy^l\\
            & \sim x^3+y^s+\widetilde{e}_0xy^k+xy^l+(\widetilde{e}_0\frac{3k-2s}{p}t^p+\widetilde{e}_1)xy^{k+p}.
        \end{aligned}
    \end{equation}
    Choosing $t$ as the solution of $\widetilde{e}_0\frac{3k-2s}{p}t^p+\widetilde{e}_1=0$, we have $$f \sim x^3+y^s+\widetilde{e}_0xy^k+xy^l\in E_{k,s,l}.$$

    \textbf{(II)} If $2s>3r$, then $in_f(f)=x^3+xy^r+y^s$. The weight vectors corresponding to the Newton diagram are $w_1=(rs,2s),\ w_2=(3rs-3r^2,3r)$ and $d=v_f(f)=3rs$. Write $f_0=in_f(f)$. Next, we find the uni-deformation system of $f$ through the calculation of $T^e_{f_0}$:

    We have $(f_0)_x=3x^2+y^{r},(f_0)_y=rxy^{r-1}+sy^{s-1}$. An easy calculation shows that the terms $\mathbf{x}^\alpha$ of the form $x^4,x^3y^k,x^2y^k,y^k$ with valuation greater than $d$ are always lied in $tj^e(f_0)$. Then we consider $xy^{r+1}$ with $v_f(xy^{r+1})=3r(s+1)$. 

    (1) If $p \nmid 2s-3r$, the equation $$\begin{pmatrix}
               1     & 3 &0 \\
            1& 1      &r  \\
          1&  0      & s \\
    \end{pmatrix}
    \begin{pmatrix}
    a\\b\\c
    \end{pmatrix}=\begin{pmatrix}0\\1\\0\end{pmatrix}$$ has a solution, which means that there exists $a,b,c$ such that $xy^{r+1}=a\cdot y\cdot f_0+b\cdot xy\cdot (f_0)_x+c \cdot y^2 \cdot (f_0)_y$. Thus, $xy^{r+1} \in tj^{e}(f_0)$ and $\Lambda=\emptyset$.

    (2) If $p \mid 2s-3r$ and $p \nmid r$, we find $$-y^{s+1-r}(f_0)_x+(\frac{3}{r}xy+y^2)(f_0)_y=xy^{r+1}(r+\frac{3s}{r}y^{s-r-1}).$$ Since $r+\frac{3s}{r}y^{s-r-1}$ is a unit of $K[[x,y]]$, we get $xy^{r+1} \in tj^{e}(f_0)$ and $\Lambda=\emptyset$.

    (3) If $p \mid 2s-3r,p \mid r$ and $s \geq 2r$, we find that $$(y+2xy)f_0-(\frac{1}{3}xy+\frac{2}{3}x^2y+y^{s-r+1})(f_0)_x=xy^{r+1}(\frac{2}{3}+2x+2y^{s-r}-\frac{2}{3}x-3xy^{s-2r}).$$ Since $\frac{2}{3}+2x+2y^{s-r}-\frac{2}{3}x-3xy^{s-2r}$ is a unit, we get $xy^{r+1} \in tj^{e}(f_0)$ and $\Lambda=\emptyset$.

    (4) If $p \mid 2s-3r,p \mid r$ and $s < 2r$, use the method of undetermined coefficients similar to Proposition \ref{example-x3ys-sudden-jump}, we can even show that there exists no $l_0,l_1,l_2 \in R$ such that $xy^{r+i}=l_0f_0+l_1(f_0)_x+l_2(f_0)_y$ for $i=1,\dots, 3r-s-1$. But we have
    $$(\frac{4}{9}+3y^{2s-3r})xy^{3r-s}=(3x+\frac{2}{3}y^{2r-s})f_0-(x^2+\frac{2}{9}xy^{2r-s}+\frac{2}{3}y^r)(f_0)_x \in tj^AC(f_0)_{6r^2}.$$
    Therefore, $\Lambda=\{xy^{r+1},\dots,xy^{3r-s-1}\}$ in this case. Using the same method as Proposition \ref{example-x3ys-SD-type}, we can show that $f$ has SD-type and separable sudden jumps $(1,k)$ with $p \mid k, r \leq k \leq 3r-s-2$.

    Using the same method as \textbf{(I)}, we can show that: if $in_f(f)=x^3+xy^r+y^s$ satisfies 
    \begin{equation}\label{x3-ys-xyr-3r-2s-4r}
        p \mid r, p \mid s, 3r<2s<4r,
    \end{equation}
    and 
    \begin{equation}\label{x3-ys-xyr-3r-2s-4r-2}
        r<r+p \leq 3r-s-2,
    \end{equation}
    then $\mathcal{K}\text{-mod}(f) \geq 2$ by Theorem \ref{main-thm}.

    If \eqref{x3-ys-xyr-3r-2s-4r} is not satisfied, using the $\alpha,\beta$-trick, we have $f \sim x^3+y^s+xy^r \in E_{0,s,r}$ by Corollary \ref{normal-cor}. Otherwise, $in_f(f)$ satisfies \eqref{x3-ys-xyr-3r-2s-4r} but does not satisfy \eqref{x3-ys-xyr-3r-2s-4r-2}, then $f \sim x^3+y^s+\widetilde{e}_0xy^r+xy^l\in E_{k,s,l}$, where $l$ satisfies $r+1 \leq l \leq 3r-s-1, p \nmid l$, $\widetilde{e}_0\in K$.

    \textbf{(III)} If $2s=3r$, we write $s=3t,r=2t$, then $in_f(f)=f_0=x^3+xy^{2t}+b_0y^{3t},\ c_0 \in K$. The weight vector corresponding to the Newton diagram $P$ is $w=(t,1)$ and $d=v_f(f)=3t$. Note that if $p\neq 31$, then $\tau^e(f_0)< \infty$. We have $$\frac{31}{6}ty^{4t-1}=\frac{2t}{3}y^{2t-1}(f_0)_x-(x-\frac{3}{2}y^t)(f_0)y $$ if $p \nmid t$. Therefore, the uni-deformation system of $f$ is given by 
    \begin{equation}
        \begin{aligned}
            &\Lambda=\{1,x,\dots,xy^{t-1},y,\dots,y^{4t-2}\}\ \\
            (resp.\ &\Lambda=\{1,x,\dots,xy^{t-1},y,\dots,y^{4t-1}\}\text{ if additionally }p \mid t).
        \end{aligned}
    \end{equation}
    when $p\neq 31$. And $\tau^e(f_0)= \infty$ if $p=31$.

    If $p \neq 31$, through calculations we know that $f$ is of SD type and has separable sudden jumps $(0,3t)$ and $(0,k)$ with $3t+1 \leq k \leq 4k-3, p \mid k-3t$. Therefore, if there exists a $k$ such that 
    \begin{equation}\label{2s=3r-mod2-condition}
        3t+1\leq k \leq 4t-3\text{ and } p\mid k-3t,
    \end{equation}
    then $G\text{-mod}(f) \geq 2$ by Theorem \ref{main-thm}.

    If \eqref{2s=3r-mod2-condition} is not satisfied, then $f \sim x^3+xy^{2t}+b_0y^{3t} \in E_{2t,3t,0}$ or $f \sim x^3+xy^{2t}+b_0y^{3t}+y^l \in E_{2t,3t,l}$ by the $\alpha,\beta$-trick, where
    \begin{equation}\label{2s=3r-condition}
        \begin{aligned}
            &t \geq 2,3t<l\leq4t-2, p \nmid l-3t \\
            &(resp.\ 3t<l\leq4t-1,p \nmid l-3t \text{ if additionally }p\mid t).
        \end{aligned}
    \end{equation}

    If $p=31$, rewrite $f$ as $f \sim x^3+xy^{2t}+b_0y^{3t}+\widetilde{b}(y)y^l$, where $\widetilde{b}(y)$ is a unit. Let $f_1=x^3+xy^{2t}+b_0y^{3t}+y^l$. Although $f_1 \neq in_f(f)=f_0$, we can replace all the content based on $f_0$ in Theorem \ref{main-thm} with that based on $f_1$, while keeping the proof process unchanged, thus obtaining a similar result. Counting separable sudden jumps, we know that if there exist a $k$ such that $3t<k \leq l-1$ and $p \mid k-3t$, then $G\text{-mod}(f) \geq 2$. Otherwise, $f \sim x^3+xy^{2t}+b'_0y^{3t}+y^l\in E_{t,l}',\ l >3t$ by $\alpha,\beta$-trick.
\end{proof}

\begin{remark}
    (i) The classification result is different from the result in fields of characteristic $0$ given by Wall. For every $p$, there are only finite many $s$ that do not satisfy \eqref{condition-k+p}, which means there are finite many case \textbf{(I)} unimodal singularities. However, Wall shows that there are infinite case \textbf{(I)} unimodal singularities in characteristic $0$ field. The main reason is: there are more sudden jumps in characteristic $p$ than in characteristic $0$, and most sudden jumps rely on $p$.\\
    (ii) The classification process is lengthy. First, we use Theorem \ref{complete-cor} to find the initial part $f_0$ of $f$. Next, we need to find a basis for the extended Tjurina algebra of $f_0$ to get the point where the jump of the extended Tjurina number occurs. Then we need to check the separability for a family of power series. After we find the bound of the modality by Theorem \ref{main-thm}, we use the implicit function theorem to finish the classification. In the following, we will omit most of the discussion and give the result directly.
\end{remark}

\begin{proposition}
    If $j_3(f) \sim x^2y$, then $f \sim x^2y+y^{s}(k\ge 4)$, which is simple.
\end{proposition}
\begin{proof}
    Write $g=x^2y$. Then $\widetilde{T}_g(\mathcal{K}g)=g+\mathfrak{m}\cdot j(g)=  \langle x^3,x^2y,xy^2 \rangle$.  Therefore, for any $l \geq 4$, we can find $$C=\mathrm{span}\langle y^j \mid  4 \leq j \leq l \rangle$$ such that $$P_{3,l} \subset C+\widetilde{T}_g(\mathcal{K}_{l}g)\cap P_{3,l}.$$ 
    By Theorem \ref{complete-cor}, we have
    $$f\sim x^2y+a(y)y^s$$
    for some $s\ge 4$, and $a(y)$ is a unit or $0$.

    If $a(y)=0$, then $f\sim x^2y$, which is not isolated. If $a(y)$ is a unit, apply the automorphism $\phi(x)=a(y)^{\frac{1}{2}}x, \phi(y)=y$. Then
    $$f\sim x^2y+a(y)y^s\sim a(y)(x^2y+y^s)\sim x^2y+y^s.$$
\end{proof}

\begin{proposition}\label{x2y+xy2}
    If $j_3(f) \sim x^2y+xy^2$, then $f \sim x^2y+y^{3}$, which is simple.
\end{proposition}
\begin{proof}
     Write $g=x^2y+xy^2$. Then $\widetilde{T}_g(\mathcal{K}g)=g+\mathfrak{m}\cdot j(g)=  \langle x^3,x^2y,xy^2,y^3 \rangle=\mathfrak{m}^3$. Therefore, $g$ is 3-determined by Theorem \ref{finite-deter} and $f\sim g= x^2y+xy^2$. Then we can apply the automorphism $\phi(x)=x+\sqrt{-1}y,\phi(y)=x-\sqrt{-1}y$.
     It follows $f\sim x^2y+y^3$.
\end{proof}

Next we discuss the case $ord(f)=4$. $j_4(f)$ is of the form $x^4,x^3y,x^2y^2,x^2y(x+y),xy(x+y)(x+ay)$ with $a\neq 0,1$.

\begin{proposition}\label{x^4}
    If $j_4(f)\sim x^4$, then $f$ belongs to the family $W$.
\end{proposition}

\begin{proof}
    Write $g=x^4$. Then $\widetilde{T}_g(\mathcal{K}g)=g+\mathfrak{m}\cdot j(g)=  \langle x^4,x^3y \rangle$. Therefore, for any $l \geq 5$, we can find $$C=\mathrm{span}\langle x^2y^{r_1},xy^{r_2},y^{r_3} \mid  3\leq r_{1} \leq l-2,4\leq r_{2} \leq l-1,5\leq r_{3}\leq l \rangle$$ such that $$P_{4,l} \subset C+\widetilde{T}_g(\mathcal{K}_{l}g)\cap P_{4,l}.$$
    By Theorem \ref{complete-cor}, we have
    $$f\sim x^4+a(y)x^2y^r+b(y)xy^s+c(y)y^t$$
    for some $r\ge 3,s\ge 4,t\ge 5$, and $a(y),b(y),c(y)$ are units or $0$. We regard $r=\infty$ (resp. $s,t=\infty$) if $a(y)=0$ (resp. $b(y),c(y)=0$).

    If $r \geq 4, s \geq 6, t \geq 8$, we write $h=j_5(f)=x^4,$ and any jet in an open neighborhood of $J_8(h)$ is of the form $h'=x^4+ax^2y^4+bxy^6+cy^8$. The codimension of $\widetilde{T}_{h'}(\mathcal{K}h') \geq 2$, which implies $\mathcal{K}\text{-mod}(f) \geq \mathcal{K}_8(f)\text{-mod}(h) \geq 2$ by Theorem \ref{cod-modality}.

    Therefore, one of the conditions $r \leq 3,s \leq 5, t \leq 7$ must be met. Note that this means $x^4+x^2y^r+xy^s+y^t$ cannot be weighted homogeneous.

    If $t=5$, then $f$ is convenient and $in_f(f)=x^4+y^5$. For $p \neq 5$, we have $f \sim x^4+y^5+\lambda x^2y^3, \lambda \in K$ by Theorem \ref{normal-cor}. Using the $\alpha,\beta$-trick, we have $f \sim x^4+y^5+\lambda x^2y^3,\ \lambda \in \{0,1\}$, which belongs to $W_{12}$ or $W_{12}'$. If $p=5$, we can show $\mathcal{K}\text{-mod}(f) \geq 2$ by Theorem \ref{cod-modality}.

    If $s=4$, we choose the $C$-polytope $P$ expanded from the Newton diagram given by $(0,4)$, $(4,1)$, $(\frac{16}{3},0)$ (the expanding point). Then $in_P(f)=x^4+xy^4$ and $f \sim x^4+xy^4+\lambda y^6,\ \lambda \in \{0,1\}$.

    If $r=3$ and $t=6$, we choose the $C$-polytope $P$ given by $(0,4),(3,2),(6,0)$ (in this case $c(y) \neq0$, otherwise $f$ is not isolated). Then $in_P(f)=x^4+x^2y^3+\lambda y^6$ and $f \sim x^4+x^2y^3+\lambda y^6+\mu y^7, \ \lambda \neq 0,\frac{1}{4},\mu\in \{0,1\}$. 

    If $r=3$ and $t \geq 7$, we choose the $C$-polytope $P$ given by $(0,4),(3,2),(t,0)$ (in this case $c(y) \neq0$, otherwise $f$ is not isolated). Then $in_P(f)=x^4+x^2y^3+\lambda y^t, \lambda \in K^{\times}$ and $f \sim x^4+x^2y^3+y^t$ (using the $\alpha,\beta$-trick, we can reduce $\lambda$).

    If $t=6$ and $r \geq 4$, then $in_P(f)=x^4+y^6$ and $f \sim x^4+y^6+\lambda x^2y^4$, where $\lambda \in \{0,1\}$.

    If $s=5$ and $p \neq 5$, we choose $P$ the expanded Newton diagram. Then $in_P(f)=x^4+xy^5$ and $f \sim x^4+xy^5+\lambda y^7+\lambda'y^8$, $\lambda,\lambda'\in K$. Using the $\alpha,\beta$-trick, we have $f \sim x^4+xy^5+\lambda y^t$, where $\lambda \in \{0,1\},k=7,8$. For case $p=5$, $\mathcal{K}\text{-mod}(f) \geq 2$ by Theorem \ref{cod-modality}.

    If $t=7$, then $in_f(f)=x^4+y^7$ and $f \sim x^4+y^7+\lambda x^2y^4+\lambda' x^2y^5$ for $p \neq 7$. Using the $\alpha,\beta$-trick, we have $f \sim x^4+y^7+\lambda x^2y^s,$ where $\lambda \in \{0,1\}$ and $s=4,5$. If $p=7$, then $\mathcal{K}\text{-mod}(f) \geq 2$ by Theorem \ref{cod-modality}. 
\end{proof}

\begin{proposition}
    If $j_4(f)\sim x^3y$, then $f$ belongs to the family $Z$.
\end{proposition}

\begin{proof}
    Write $g=x^3y$. Then $\widetilde{T}_g(\mathcal{K}g)=g+\mathfrak{m}\cdot j(g)=  \langle x^4,x^3y,x^2y^2 \rangle$. Therefore, for any $l \geq 5$, we can find $$C=\mathrm{span}\langle xy^i,y^j \mid  4 \leq i \leq l-1,5 \leq j \leq l \rangle$$ such that $$P_{4,l} \subset C+\widetilde{T}_g(\mathcal{K}_{l}g)\cap P_{4,l}.$$
    By Theorem \ref{complete-cor}, we have
    $$f\sim x^3y+a(y)xy^r+b(y)y^s$$
    for some $r\ge 4,s\ge 5$, and $a(y),b(y)$ are units or $0$.
    
    If $a(y)=0$, $f \sim x^3y+b(y)y^s,\ s \geq 5$. Applying the automorphism $\phi(x)=b(y)^{\frac{1}{3}}x, \phi(y)=y$, we have $$f \sim b(y)(x^3y+y^s) \sim x^3y+y^s,\ s \geq 5.$$

    If $b(y)=0$, similarly we have $f \sim x^3y+xy^r,\ r \geq 4$.

    Next, we assume that $a(y),b(y)$ are both units. Then the Newton diagram and $in_P(f)$ depend on $r,s$. This case is similar to Proposition \ref{x^3}.

    \textbf{(I)} If $2s+1<3r$, we expand the Newton diagram to get the $C$-polytope $P$, which is given by $(0,\frac{3s}{s-1})$ (the expanding point)$,(1,3),(s,0)$. Then $in_P(f)=x^3y+y^s$. The weight vector corresponding to $P$ is $(s-1,3)$ and $d=v_P(f)=3s$.

    Similarly to case \textbf{(I)} in Proposition \ref{x^3}, we have the following:\\
    If there exists a $k$ such that $$\left\lfloor {\frac{2s+1}{3}} \right\rfloor+1 \leq k \leq k+p \leq s-3$$ and $p\mid 3k-2s-1$, then $\mathcal{K}\text{-mod}(f)\geq 2$. Otherwise, we have $$f \sim x^3y+y^s+xy^k,\ p \nmid 3k-2s-1$$ or $$f\sim x^3y+y^s+\widetilde{e}_0xy^k+xy^l,\ l>k,p \mid 3k-2s-1, p \nmid 3l-2s-1.$$

    \textbf{(II)} If $2s+1>3r$, we expand the Newton diagram to get the $C$-polytope $P$, which is given by $(0,\frac{3r-1}{r-1})$(the expanding point)$,(1,3),(r,1),(s,0)$. Then $in_P(f)=x^3y+xy^r+y^s$. The weight vectors corresponding to $P$ are $w_1=((r-1)s,2s),\ w_2=((3r-1)(s-r),3r-1)$ and $d=v_P(f)=(3r-1)s$. 

    Similarly to case \textbf{(II)} in Proposition \ref{x^3}, we have the following:\\
    If $p \mid 3r-2s-1$ and $r<r+p \leq 3r-s$, then $\mathcal{K}\text{-mod}(f) \geq 2$. Otherwise, $f \sim x^3y+y^s+xy^r$ for $p \nmid 3r-2s-1$ or $f \sim x^3y+y^s+\widetilde{e}_0xy^r+xy^l$, where $p \mid 3r-2s-1,r+1 \leq l \leq 3r-s+1,p \nmid 3l-s-1$.

    \textbf{(III)} If $2s+1=3r$, we write $s=3t+1,r=2t+1$ and expand the Newton diagram to get the $C$-polytope $P$, which is given by $(0,\frac{3t+1}{t})$(the expanding point)$,(1,3),(2t+1,1),(3t+1,0)$. Then $in_P(f)=f_0=x^3y+xy^{2t+1}+b_0y^{3t+1}$.

    Similarly to case \textbf{(III)} in Proposition \ref{x^3}, we have the following:\\
    For $p \neq 31$, if there exists a $k$ such that $3t+2 \leq k\leq 4t-1$ and $p \mid k-3t-1$, then $\mathcal{K}\text{-mod}(f) \geq 2$. Otherwise, $f \sim x^3y+xy^{2t+1}+b_0y^{3t+1}+y^l$, where
    \begin{equation}
        t \geq 2, 3t+2 \leq l \leq 4t,p \nmid l-3t-1.
    \end{equation}
    For $p=31$, we have $f \sim x^3y+xy^{2t+1}+b_0y^{3t+1}+y^l$, where $l>3t+1$ and there does not exist $k$ such that $3t+2\leq k \leq l$ and $p \nmid k-3t-1$. 
\end{proof}

\begin{proposition}
    If $j_4(f)\sim x^2y^2$, then $f\sim x^s+x^2y^2+y^t \in T_{r,s,2},\ s,t\geq 5$.
\end{proposition}

\begin{proof}
    Write $g=x^2y^2$. Then $\widetilde{T}_g(\mathcal{K}g)=g+\mathfrak{m}\cdot j(g)=  \langle x^3y,x^2y^2,xy^3 \rangle$. Therefore, for any $l \geq 5$, we can find $$C=\mathrm{span}\langle x^i,y^j \mid  5 \leq i \leq l,5 \leq j \leq l \rangle$$ such that $$P_{4,l} \subset C+\widetilde{T}_g(\mathcal{K}_{l}g)\cap P_{4,l}.$$ 
    By Theorem \ref{complete-cor}, we have
    $$f\sim x^2y^2+a(x)x^s+b(y)y^t$$
    for some $s,t\ge 5$, and $a(x),b(y)$ is a unit or $0$. In reality, $f$ is not isolated if $a(x)=0$ or $b(y)=0$. 
    
    Therefore, $in_f(f)=x^s+x^2y^2+y^t$ and $d=v_f(f)=2rs$. Calculation shows that the regular basis $B$ of $T_f$
    is contained in $\{1,x,x^2,\ldots,x^s,y,y^2,\ldots, y^t,xy\}$ for all $s,t \geq 5$. Since $v_P({\bf x}^\alpha)\le d$ for all ${\bf x}^\alpha\in B$, $f\sim x^s+x^2y^2+y^t$ by Corollary \ref{normal-cor}.
\end{proof}

\begin{proposition}
    If $j_4(f)\sim x^2y(x+y)$, then $f\sim x^4+x^2y^2+y^s \in T_{4,s,2}$, $s \geq 5$.
\end{proposition}
\begin{proof}
    We have $j_4(f)\sim x^2y(x+y)\sim x^2(x^2+y^2). $ Write $g=x^4+x^2y^2$. Then $\widetilde{T}_g(\mathcal{K}g)=g+\mathfrak{m}\cdot j(g)=  \langle x^4,x^3y,x^2y^2,xy^3 \rangle$. Therefore, for any $l \geq 5$, we can find $$C=\mathrm{span}\langle y^j \mid  5 \leq j \leq l \rangle$$ such that $$P_{4,l} \subset C+\widetilde{T}_g(\mathcal{K}_{l}g)\cap P_{4,l}.$$ 
    By Theorem \ref{complete-cor}, we have
    $$f\sim x^4+x^2y^2+a(y)y^s$$
    for some $s\ge 5$,and $a(y)$ is a unit or $0$. $f$ is not isolated when $a(y)y^s=0$, so $a(y)\neq 0$. Therefore $in_f(f)=x^4+x^2y^2+y^s$. By Corollary \ref{normal-cor}, we have
    $f\sim x^4+x^2y^2+y^s,\ s\geq 5$.    
\end{proof}

\begin{proposition}
    If $j_4(f)\sim xy(x+y)(x+ay)$ with $a\neq 0,1$, then $f\sim x^4+y^4+bx^2y^2\quad(b^2\neq 4) \in T_{4,4,2}$.
\end{proposition}
\begin{proof}
    Write $g=xy(x+y)(x+ay)$. We now show $$\widetilde{T}_g(\mathcal{K}g)= g+\mathfrak{m}\cdot j(g)\supset  \langle x^5,x^4y,x^3y^2,x^2y^3,xy^4,y^5 \rangle=\mathfrak{m}^5$$
    through the following calculations. Note that 
    $$
    \begin{pmatrix}y^2g_x\\xyg_x\\x^2g_x\\y(2g-xg_x)\\x(2g-xg_x)\end{pmatrix}
    =\begin{pmatrix}
       &        & 3     & 2(a+1) & a\\
       & 3      & 2(a+1)& a      &  \\
    3  & 2(a+1) & a     &        &  \\
       & -1     &       & a      &  \\
    -1 &        & a     &        &  \\
    \end{pmatrix}
    \begin{pmatrix}
    x^4y\\x^3y^2\\x^2y^3\\xy^4\\y^5
    \end{pmatrix},
    $$
    where the determinant of the coefficient matrix on the right-hand side is equal to $-4a^2(a-1)^2\neq0$,\ so we have $x^4y,x^3y^2,x^2y^3,xy^4,y^5\in \widetilde{T}_g(\mathcal{K}g)=g+\mathfrak{m}\cdot j(g)$. Furthermore, the identity $x^5=x^2g_{y}-2(a+1)x^4y-3ax^3y^2$ implies that $x^5\in \widetilde{T}_g(\mathcal{K}g)=g+\mathfrak{m}\cdot j(g)$, which shows $\widetilde{T}_g(\mathcal{K}g)\supset\mathfrak{m}^5$. By Theorem \ref{complete-cor}, we have
    $$f\sim xy(x+y)(x+ay).$$

    This expression will be transformed into a canonical representation through the calculations below. Let $\lambda$ be the root of the equation $\lambda^2+2\lambda+\frac{1}{a}=0$ and let $t^2=-\frac{\lambda}{\lambda+2}$. Then
    $$
    \begin{aligned}
    f&\sim xy(x+y)(x+ay) \qquad \left( x\mapsto \frac{x}{\lambda},y\mapsto y \right) \\
    &\sim \left( \frac{x}{\lambda} \right)y\left( \frac{x}{\lambda}+y \right)\left( \frac{x}{\lambda}+ay \right) \\
    &\sim xy(x+\lambda y)(x+\lambda ay) \qquad \left( x\mapsto x+ty,y\mapsto x+\frac{y}{t}\right) \\
    &\sim (x+ty)\left( x+\frac{y}{t} \right)\left( (1+\lambda)x+\left( t+\frac{\lambda}{t} \right)y \right)\left( (1+\lambda a)x+\left( t+\frac{\lambda a}{t} \right)y \right) \\
    &\sim (x+ty)\left( x+\frac{y}{t} \right)\left( x+ \frac{t^2+\lambda}{t(1+\lambda)} y \right)\left(x+\frac{t^2+\lambda a}{t(1+\lambda a)} y \right).
    \end{aligned}
    $$
    Since $\lambda=-\frac{-2t^2}{t^2+1}$, we have
    $$\frac{t^2+\lambda}{t(1+\lambda)}=\frac{t^2-\frac{-2t^2}{t^2+1}}{t(1-\frac{-2t^2}{t^2+1})}=-t,$$
    $$\frac{t^2+\lambda a}{t(1+\lambda a)} =\frac{t^2-\frac{1}{\lambda+2}}{t(1-\frac{1}{\lambda+2})}=\frac{(\lambda+2)t^2-1}{t(\lambda+1)} =\frac{-\lambda-1}{t(\lambda+1)}=-\frac{1}{t},$$
    so
    $$
    \begin{aligned}
    f&\sim (x+ty)\left( x+\frac{y}{t} \right)(x-ty)\left( x-\frac{y}{t} \right) \\
    &=x^4+y^4-\left( t^2+\frac{1}{t^2} \right)x^2y^2.
    \end{aligned}
    $$
    Denote $-t^2-\frac{1}{t^2} $ by $b$. Then 
    $$
    \begin{aligned}
    b&=\frac{\lambda}{\lambda+2}+\frac{\lambda+2}{\lambda} =\frac{2(\lambda^2+2\lambda)+4}{\lambda^2+2\lambda} =\frac{2\left( -\frac{1}{a} \right)+4}{-\frac{1}{a}} =2-4a\in K\backslash\{2,-2\}
    \end{aligned}$$
    and
    $$f\sim x^4+y^4+bx^2y^2\quad(b^2\neq 4).$$
\end{proof}

\subsection{Unimodal hypersurface singularities in $K[[x_1,\dots,x_n]]$ with order 2}
For $f \in \mathfrak{m}^2 \subset K[[x_1,\dots,x_n]]$, assume $n \geq 3,\ l=\mathrm{ord}(f)=2$.

By  Lemma \ref{splitting}, we have $f({\bf x}) \sim x_1^2+g({\bf x'}),$ where ${\bf x'}=(x_2,\dots,x_n)$. In fact, we have:
\begin{lemma}\label{split-unique}
    Let $f_1({\bf x})=x_1^2+g_1({\bf x'})$, $f_2({\bf x})=x_1^2+g_2({\bf x'})$. Then $f_1 \sim f_2 \Longleftrightarrow g_1 \sim g_2$.    
\end{lemma}

To prove Lemma \ref{split-unique}, we need the Mather-Yau Theorem in positive characteristic:
\begin{definition}
    Define $T_k(f)=K[[{\bf x}]]/\langle f, \mathfrak{m}^k \cdot j(f) \rangle$ as the $k$-th Tjurina algebra, where $j(f)=\langle \frac{\partial f}{\partial x_1}, \dots, \frac{\partial f}{\partial x_n} \rangle$ is the Jacobi ideal.
\end{definition}

\begin{theorem}[\cite{mather-yau} Theorem 2.2]\label{mather-yau}
Let $f,g\in K [[\bf{x} 
 ]]$ be such that $ord(f)=s \geq 2$ and $\tau(f)<\infty$. 
Then the following are equivalent:\\
i) $f \sim g$. \\
ii) $T_k(f)  \cong T_k(g)$ as $K$-algebras for some (equivalently for all) $k$ such that $$\mathfrak{m}^{\left\lfloor {\frac{{k + 2s}}{2}} \right\rfloor}\subset\mathfrak{m}\cdot \widetilde{T}_f(\mathcal{K}f)$$
\noindent where ${\left\lfloor {\frac{{k + 2s}}{2}} \right\rfloor}$ means the maximal integer which does not exceed $\frac{{k + 2s}}{2}$.
\end{theorem}

Then we can begin the proof of Lemma \ref{split-unique}.
\begin{proof}
    First, assume $g_1 \sim g_2$. Then there exist $u({\bf x'}) \in K[[\bf x']]^{\times}$ and $\Phi' \in \mathrm{Aut}K[[\bf x']]$ such that $g_2({\bf x'})=u({\bf x'}) \cdot g_1(\Phi'({\bf x'}))$. Then we can apply $$\Phi \in \mathrm{Aut}K[[{\bf x}]]:\ x_1 \mapsto u({\bf x'})^{-\frac{1}{2}}x_1,\ {\bf x'} \mapsto \Phi'({\bf x'}).$$ It follows that
    $$u({\bf x'})\Phi(f_1)=u({\bf x'})f_1(\Phi({\bf x}))=x_1^2+u({\bf x'}) \cdot g_1(\Phi'({\bf x'}))=x_1^2+g_2({\bf x'})=f_2. $$  This implies $f_1 \sim f_2$.

    Next, we assume $f_1 \sim f_2$. By Theorem \ref{mather-yau}, there exists a $k \in \mathbb{N}$ such that
    \begin{equation}\label{deter-condition}
        \mathfrak{m}^{\left\lfloor {\frac{{k + 4}}{2}} \right\rfloor}\subset\mathfrak{m}\cdot \widetilde{T}_{f_1}(\mathcal{K}{f_1})
    \end{equation}
    and
    \begin{equation}\label{algebra-condition}
        T_k(f_1)  \cong T_k(f_2),
    \end{equation}
    i.e. $$K[[{\bf x}]]/\langle x_1^2+g_1, \mathfrak{m}^k \cdot x_1, \mathfrak{m}^k \cdot j(g_1) \rangle \cong K[[{\bf x}]]/\langle x_1^2+g_2,\mathfrak{m}^k\cdot x_1, \mathfrak{m}^k \cdot j(g_2) \rangle .$$
    
    Modulo $\langle x_1\rangle$ on both sides of \eqref{algebra-condition} and write $\mathfrak{m}'=\langle x_2,\dots,x_n\rangle \subset K[[{\bf x'}]]$, we have
    $$T_k(g_1) \cong K[[{\bf x'}]]/\langle g_1, (\mathfrak{m}')^k \cdot j(g_1) \rangle \cong K[[{\bf x'}]]/g_2, (\mathfrak{m}')^k \cdot j(g_2) \rangle \cong T_k(g_2) .$$
    Similarly, modulo $\langle x_1\rangle$ on both sides of \ref{deter-condition}, we have $$\mathfrak{m}'^{\left\lfloor {\frac{{k + 2\mathrm{ord}(g_1)}}{2}} \right\rfloor}\subset \mathfrak{m}'^{\left\lfloor {\frac{{k + 4}}{2}} \right\rfloor}\subset\mathfrak{m}'\cdot \widetilde{T}_{g_1}(\mathcal{K}{g_1}).$$
    By Theorem \ref{mather-yau} again, we get $g_1 \sim g_2$.    
\end{proof}

\begin{corollary}\label{split-unique-cor}
    $g(x_1,\ldots, x_k)$ in Lemma \ref{splitting} is unique up to contact equivalence.
\end{corollary}

Therefore, we can show
\begin{proposition}\label{split-modality}
    The $\mathcal{K}$\text{-modality} of $f({\bf x})$ in $K[[{\bf x}]]$ is equal to the $\mathcal{K}$\text{-modality} of $g({\bf x'})$ in $K[[{\bf x'}]]$.
\end{proposition}
\begin{proof}
    Using the same argument as in the proof \cite[Lemma 3.11]{right-simple}.
\end{proof}

Combining Corollary \ref{split-unique-cor} and Proposition \ref{split-modality}, we need only to consider the classification of unimodal singularities $g(x_1,\dots, x_k) \in K[[x_1,\dots, x_k]]$ with $k<n$ and $\mathrm{ord}(g) \geq 3$. Moreover, as a result of Proposition \ref{nl-bound}, we can easily prove Proposition \ref{class-n>3}. Thus, we only need to classify the unimodal isolated hypersurface singularity with $n=3,l=3$. 

\subsection{Unimodal hypersurface singularities in $K[[x,y,z]]$ with order 3}
As shown in \cite{right_unimodal}, the $3$-jets in $K[[x,y,z]]$ are contact equivalent to the following form:
\begin{equation}\label{xyz-3-jet}
\begin{aligned}
    &x^3+y^3+z^3+axyz\ (a^3+27 \neq 0),\ x^3+y^3+xyz,\ x^3+xyz,\ xyz,\\
    &x^3+yz^2,\ x^2z+yz^2,\ x^3+xz^2,\ x^2y,\ x^3.
\end{aligned}
\end{equation}

One can show that $x^3+y^3+z^3+axyz\ (a^3+27 \neq 0)$ is $3$-determined, therefore the corresponding normal form is $x^3+y^3+z^3+axyz\ (a^3+27 \neq 0)$.

\begin{proposition}
    If $j_3(f)$ is of the form $x^3+y^3+xyz,\ x^3+xyz,\ xyz,$ then $f$ is contact equivalent to $x^r+y^s+z^t+xyz$ for $r,s \geq 3$, $t \geq 4$. That is, $f$ belongs to the family $T_{r,s,t}$. 
\end{proposition}
\begin{proof}
    If $j_3(f) \sim x^3+y^3+xyz$ (resp. $x^3+xyz,xyz$), the complete transversal is given by $$C=\text{span}\langle z^4,z^5,\dots \rangle$$ $$(resp.\ \text{span}\langle y^4,y^5,\dots,z^4,z^5,\dots \rangle,\ \text{span}\langle x^4,\dots,y^4,\dots,z^4,\dots \rangle).$$ Therefore, $f \sim a(x)x^r+b(y)y^s+c(z)z^t+xyz$, where $r,s \geq 3$, $t \geq 4$. Note that if one of $a(x),b(y),c(z)$ is $0$, then $f$ is not isolated. Hence $a(x),b(y),c(z)$ are all units and $in_f(f)=x^r+y^s+z^t+xyz,\ d=v_f(f)=rst$. In addition, there are no terms in a basis of $gr^{AC}_P(T_{in_f(f)})$ with valuation greater than $d$. By Corollary \ref{normal-cor}, we have $f \sim x^r+y^s+z^t+xyz$.
\end{proof}

\begin{proposition}
    If $j_3(f) \sim x^3+yz^2$, then $f$ belongs to the family $Q$.
\end{proposition}

\begin{proof}
    By Theorem \ref{complete-cor}, we have $f \sim x^3+yz^2+a(y)xy^r+b(y)y^s,\ r \geq3,s \geq 4$.

    If $a(y)=0$, then $f \sim x^3+yz^2+b(y)y^s \sim x^3+yz^2+y^s,\ s \geq 4$ using the $\alpha,\beta$-trick. If $b(y)=0$, similarly $f \sim x^3+yz^2+xy^r,\ r \geq 3$.

    Next we assume that $a(y),b(y)$ are all units. The Newton diagram depends on $r,s$.
    
    \textbf{(I)} If $2s<3r$, we choose the $C$-polytope expanding from the Newton diagram as $(3,0,0)$, $(0,s,0)$, $(0,0,3)$(the expanding point), $(0,1,2)$. Then $in_P(f)=x^3+y^s+yz^2$. This case is very similar to the case \textbf{(I)} in Proposition \ref{x^3}. The uni-deformation system of $f$ is given by $\Lambda=\{(1,\left\lfloor {\frac{2}{3}s} \right\rfloor+1,0),\dots,(1,s-2,0)\}$ for $3 \nmid s$ and $\Lambda=\{(1,\frac{2}{3}s,0),\dots, (1,s-2,0)\}$ for $3 \mid s$ (and additionally $(1,{s-1},0)$ if $p \mid s$). And we have the same result: if there exists a $k$ such that $$\left\lfloor {\frac{2}{3}s} \right\rfloor+1 \leq k \leq k+p \leq s-3$$ and $p\mid 3k-2s$, then $\mathcal{K}\text{-mod}(f)\geq 2$. Otherwise, $f \sim x^3+y^s+yz^2+\widetilde{e}_0xy^k+xy^l$, where $s \geq 4, p \mid 3k-2s, l>k,p \nmid 3l-2s$ and $\widetilde{e}_0 \in K$.

    \textbf{(II)} The case $2s>3r$. This is similar to Proposition \ref{x^3}, case \textbf{(II)}. If $r,s$ satisfies \eqref{x3-ys-xyr-3r-2s-4r} and \eqref{x3-ys-xyr-3r-2s-4r-2}, then $\mathcal{K}\text{-mod}(f)\geq 2$. Otherwise, $f \sim x^3+y^s+yz^2+\widetilde{e}_0xy^r+xy^l$, where $r+1\leq l \leq 3r-s-1, p \nmid 3l-2s$ and $\widetilde{e}_0 \in K$.

    \textbf{(III)} The case $2s=3r$. This is similar to Proposition \ref{x^3}, case \textbf{(III)}. Write $s=3t,r=2t$. For $p \neq 31$, if \eqref{2s=3r-mod2-condition} is satisfied, then $\mathcal{K}\text{-mod}(f)\geq 2$. Otherwise, $f \sim x^3+xy^{2t}+b_0y^{3t}+y^l+yz^2$, where $l,t$ satisfies \eqref{2s=3r-condition}. For $p=31$, $f \sim x^3+xy^{2t}+b_0y^{3t}+y^l+yz^2,\ l>3t$ such that there does not exist a $k$ with $3t<k<l$ and $p \mid k-3t$.  
\end{proof}

\begin{proposition}\label{x2z+yz2}
    If $j_3(f) \sim x^2z+yz^2$, then $f$ belongs to the family $S$.
\end{proposition}
\begin{proof}
    A complete transversal $C$ is given by $\{x^2y^2,xy^3,y^4,x^2y^3,\dots\}$. Theorem \ref{complete-cor} shows $f \sim x^2z+yz^2+a(y)x^2y^r+b(y)xy^s+c(y)y^t$, where $r \geq 2,s \geq 3,t \geq 4$.

    This case is similar to Proposition \ref{x^4}. If $r \geq 3, s \geq 5, t\geq 7$, we have $\mathcal{K}\text{-mod}(f) \geq 2$.

    For the rest of the cases, we have the following:

    If $t=4$, $f \sim x^2z+yz^2+y^4+\lambda x^2y^2,\ \lambda \in \{0,1\}$.

    If $s=3$, $f \sim x^2z+yz^2+xy^3$.

    If $r=2$ and $s=4$, $f \sim x^2z+yz^2+x^2y^2+xy^4$.

    If $r=2$ and $t= 5$, $f \sim x^2z+yz^2+x^2y^2+\lambda y^5+\mu y^6,\ \lambda \neq 0,-1, \mu \in \{0,1\}$.

    If $r=2$, $s\geq 5,s+2\leq t\leq 2s-3$, $f \sim x^2z+yz^2+x^2y^2+xy^s+\lambda y^t,\ \lambda \neq 0$.

    If $r=2,\ s\geq 5,\ t>2s-3$, $f \sim x^2z+yz^2+x^2y^2+xy^s$.

    If $r=2,\ 6\leq t<s+2$, $f \sim x^2z+yz^2+x^2y^2+y^t$.

    If $r \geq 3,s=4$, $f \sim x^2z+yz^2+xy^4+\lambda y^t,\ \lambda \in \{0,1\}, t=6,7$.

    If $r \geq 3, s\geq 5,t=5$, $f \sim x^2z+yz^2+\lambda x^2y^3+y^5, \lambda \in \{0,1\}$ for $p \neq 5$. If $p=5$, then $\mathcal{K}\text{-mod}(f) \geq 2$.

    If $r \geq 3, s\geq 5,t=6$, $f \sim x^2z+yz^2+\lambda x^2y^k+y^6, \lambda \in \{0,1\},k=3,4$.
\end{proof}

\begin{proposition}
    If $j_3(f) \sim x^3+xz^2$, then $f$ belongs to the family $U$.
\end{proposition}
\begin{proof}
    Note that $f \sim x^3+xz^2+a(y)x^2y^r+b(y)xy^s+c(y)y^tz+d(y)y^w$ for $r \geq 2,s \geq 3,t \geq 3,w \geq 4$ and $\mathcal{K}\text{-mod}(f) \geq 2$ if $r \geq 2,s \geq 4,t\geq 4,w\geq 6$. By a similar discussion to Proposition \ref{x2z+yz2}. We get $f$ is contact equivalent to the following forms:
    $$x^3+xz^2+xy^3+y^tz,\ t\geq 4;\ x^3+xz^2+y^4+\lambda x^2y^2,\ \lambda \in  \{0,1\}; $$
    $$x^3+xz^2+y^5+\lambda x^2y^3,\ \lambda \in  \{0,1\} \text{ for }p \neq 5;$$
    $$x^3+xz^2+xy^3+\lambda y^3z+\mu y^4z,\ \lambda^2\neq 0,-1,\mu \in \{0,1\};$$
    $$x^3+xz^2+y^3z+\lambda xy^4,\ \lambda \in \{0,1\}.$$
\end{proof}

For the last two cases in \eqref{xyz-3-jet}, using Theorem \ref{cod-modality} we can show:
\begin{proposition}
    If $j_3(f) \sim x^2y$ or $x^3$, then $\mathcal{K}\text{-mod}(f) \geq 2$.
\end{proposition}

So far, we have finished the proof of Proposition \ref{class-xy} to Proposition \ref{class-n>3}.

\section{Check the modality}
In this section, we will check whether the candidates in Table \ref{class-xy} and Table \ref{table-xyz} are unimodal. We have the following propositions of the modality from \cite{right-simple}:

\begin{proposition}\label{check-mod}
    For $f \in \mathfrak{m}$ being a power series such that $\tau(f)<\infty$. Let $$F(\mathbf{t},\mathbf{x})=f(\mathbf{x})+\sum_{i=1}^d t_ig_i(\mathbf{x}),$$ where $g_i$ is a $K$-basis of $T^{e,sec}_f=\mathfrak{m}/\widetilde{T}_f(\mathcal{K}f)$ and $\mathbf{t}=(t_1,\dots,t_d) \in T=\mathrm{Spec}K[[t_1,\dots,t_d]]$. $F(\mathbf{t},\mathbf{x})$ is called the semiuniversal deformation of $f$.\\
    (1) By a $\mathcal{K}$-modular family over a subvariety $S$ of $T=\mathrm{Spec}K[[t_1,\dots,t_d]]$, we mean a family $h_s(x) \in \mathcal{O}(S)[[\mathbf{x}]]$ such that for every $s \in S$, there is only finitely many $s' \in S$ such that $h_{s'}\sim h_s$.\\
    (2)Assume that there exist an open neighborhood $W\subset T$ of $0$ and $\mathcal{K}$-modular families $h_{s_i}^{(i)}(\mathbf{x})$, $i=1,\dots,q$ and that for each open neighborhood $V \subset W$ of $0$ and for all $s_i \in S_i$ there exist a $\mathbf{t} \in V$ such that $F(\mathbf{x},\mathbf{t}) \sim h_{s_i}^{(i)}(\mathbf{x})$, then $\mathcal{K}\text{-mod}(f) =\max_{i=1,\dots,q} \{\dim S_i\}$.
\end{proposition}

\begin{proposition}\label{semi-cont-mod}
    The $\mathcal{K}$\text{-modality} is upper semicontinunous. That is, for all $i \in \mathbb{N}$, the sets $$U_i=\{f \in \mathfrak{m} \subset K[[\mathbf{x}]]|\mathcal{K}\text{-mod}(f)\leq i\}$$ are open in $K[[\mathbf{x}]]$. Moreover, for $f,T$ and $F(\mathbf{t},\mathbf{x})$ defined above, the set $$\{\mathbf{t}\in T|\mathcal{K}\text{-mod}(F(\mathbf{t},\mathbf{x}))\leq \mathcal{K}\text{-mod}(f)\} $$ is open in $T$.
\end{proposition}

By Proposition \ref{check-mod} and Proposition \ref{semi-cont-mod}, we only need to consider the semiuniversal deformation of normal forms in Table \ref{class-xy} and Table \ref{table-xyz} and show that they can only deform to families with dimensions of $0$ or $1$. We calculate the family $E$, for example.

For type $E_{0,s}$ of the form $f=x^3+y^s$, a basis of $T^{e,sec}_f$ is given by $$\{x,x^2,xy\dots,xy^{s-2},y,\dots,y^{s-1}\}$$ (resp. $\{x,x^2,xy\dots,xy^{s-1},y,\dots,y^{s-1}\}$ if $p \mid s$). 

If $p \nmid s$, we have $$F(\mathbf{x},\mathbf{t})=x^3+y^s+t_1x+t_2x^2+t_3xy+\dots+t_{s}xy^{s-2}+t_{s+1}y+\dots+t_{2s-1}y^{s-1}.$$ 
(1) If $t_1,t_{s+1} \neq 0$, then $F(\mathbf{x},\mathbf{t})$ is not singular.\\
(2) If $t_1,t_{s+1}=0$ and $t_2 \neq 0$, by Lemma \ref{splitting}, we have $F(\mathbf{x},\mathbf{t}) \sim x^2+g(y,\mathbf{t})$, which is simple (of modality 0). Similarly $t_3,t_{s+2}=0$.\\
(3) If $t_1,t_2,t_3,t_{s+1},t_{s+2}=0$ and $t_4 \neq 0$, then $j_3(F(\mathbf{x},\mathbf{t})) \sim x^3+t_4xy^2\sim x^2y+xy^2$, which is simple by Proposition \ref{x2y+xy2}. Similarly $t_{s+3}=0$.\\
(4) If $t_1,\dots,t_4=0,t_{s+1},\dots, t_{s+3}=0$ and $t_{s+4} \neq 0$, we denote $g=F(\mathbf{x},\mathbf{t})$. Then $in_g(g) \sim x^3+t_{s+4}y^4$ and $g \sim x^3+t_{s+4}y^4\sim x^3+y^4$ by Corollary \ref{normal-cor}, which is simple. Similarly, $t_{s+5},t_{5}=0.$\\
(5) If $t_1,\dots,t_5=0,t_{s+1},\dots, t_{s+5}=0$ and $t_{s+4} \neq 0$ and $t_{6} \neq 0$, we denote $g=F(\mathbf{x},\mathbf{t})$. Moreover, assume $t_{s+6} \neq 0$. Then $in_g(g)=x^3+t_6xy^4+t_{s+6}y^6 \sim x^3+xy^4+\lambda y^6$, $\lambda \neq 0$. If $p \neq 31$, then $g \sim x^3+xy^4+\lambda y^6$ by Corollary \ref{normal-cor}, which is a family of $\dim 1$. For the case $p=31$, if there exists a $k$ such that $6<k\leq s-2$ and $p \mid k-6$, then $\mathcal{K}\text{-mod}(x^3+xy^4+\lambda y^6+t_{2s-1}y^{s-1}) \geq 2$, which means that $\mathcal{K}\text{-mod}(f) \geq 2$ by Proposition \ref{x^3}, case \textbf{(III)}. If such $k$ does not exist, then $g \sim x^3+xy^4+\lambda y^6+y^l$ for some $6<l \leq s-1$ by Proposition \ref{x^3}, case \textbf{(III)}, which is a family of $\dim 1$. Hence $f$ is simple.\\
(6) Similarly, if there exist $u,v$ such that $3 \leq u<v \leq s-1$ and one of the following
\begin{equation}
    \begin{aligned}
        &A:\ \left\lfloor {\frac{2}{3}v} \right\rfloor+1 \leq u \leq u+p \leq v-3,\ p\mid 3u-2v;\\
        &B:\ p \mid u, p \mid v, 3u<2v<4u, u<u+p\leq 3u-v-2;\\
        &C:\ p \neq 31, u \mathrm{\ is\ even},\frac{3}{2}u+1\leq v\leq 2u-3,p \mid v-\frac{3}{2}u;\\
        &D:\ p=31,u \mathrm{\ is\ even},\frac{3}{2}u+1\leq v\leq s-2,p \mid v-\frac{3}{2}u;
    \end{aligned}
\end{equation}
holds, then $\mathcal{K}\text{-mod}(f) \geq 2$. Otherwise, $\mathcal{K}\text{-mod}(f)\leq 1$.

Using this method, we can present all types of unimodal hypersurface singularities.

\begin{theorem}\label{final-class}
    Let $K$ be an algebraically closed field of characteristic $p>3$. Then every unimodal hypersurface singularity is contact equivalent to one of the following forms:
\renewcommand\arraystretch{1.5}
\begin{longtable}{|c|c|c|}
\caption{}\label{final-table-xy}
\\
\hline
Symbol&Form& condition  \\
\hline
$E_{0,s}$&$x^3+y^{s}$&$s \geq 6$ and do not exist $3\leq u<v\leq s-1$\\
&&(resp. $3\leq u<v\leq s$ if additionally $p \mid s$) \\
&&such that any of the condition \eqref{condition-x3} holds \\
&&\\
\hline
$E_{r,0}$&$x^3+xy^{r}$&$r \geq 4$ and do not exist $3\leq u\leq r-1,4\leq v\leq 2r-2$ \\
&& (resp. $4\leq v\leq 2r-1$ if additionally $p \mid r$)\\
&&such that any of the condition \eqref{condition-x3} holds \\
\hline
$E_{r,s}^0$&$x^3+y^s+xy^r$& $s\geq 4,\frac{2}{3}s<r \leq s-2,p\nmid 3r-2s$\\
&& (resp. $\frac{2}{3}s<r \leq s-1$ if additionally $p \mid s$)\\
&& and do not exist $3\leq u\leq r-1,4\leq v \leq s-1$\\
&&such that any of the condition \eqref{condition-x3} holds \\
\hline
$E_{r,s}^{0'}$&$x^3+y^s+xy^r$& $s\geq 4,\frac{2}{3}s<r \leq s-2,p\mid 3r-2s$\\
&except for $x^3+xy^4+y^5$& (resp. $\frac{2}{3}s<r \leq s-1$ if additionally $p \mid s$)\\
&when $p=5$& and do not exist $3\leq u\leq s-2,4\leq v \leq s-1$\\
&(which is simple)& (resp. $3\leq u\leq s-2$ if additionally $p \mid s$)\\
&&such that any of the condition \eqref{condition-x3} holds \\
\hline
$E_{r,s}^1$&$x^3+y^s+xy^r$& $r\geq 3,3r<2s<4r,p\nmid 3r-2s$\\
&& and do not exist $3\leq u\leq r-1,4\leq v \leq s-1$\\
&&such that any of the condition \eqref{condition-x3} holds \\
\hline
$E_{r,s}^{1'}$&$x^3+y^s+xy^r$& $r\geq 3,3r<2s<4r,p\mid 3r-2s$\\
&& and do not exist $3\leq u\leq r-1,4\leq v \leq s$\\
&&(resp. $3\leq u \leq 3r-s-1$ if additionally $p \mid r,s$)\\
&&such that any of the condition \eqref{condition-x3} holds \\
\hline
$E_{k,s,l}^0$&$x^3+y^s+\lambda xy^k+xy^l$& $s\geq 4,\frac{2}{3}s<k<l\leq s-2,p\mid 3k-2s,p \nmid 3l-2s,\lambda \neq 0$\\
&&(resp. $\frac{2}{3}s<k<l\leq s-1$ if additionally $p \mid s$)\\
&& and do not exist $3\leq u\leq l-1,4\leq v \leq s-1$\\
&&such that any of the condition \eqref{condition-x3} holds \\
\hline
$E_{k,s,l}^1$&$x^3+y^s+\lambda xy^k+xy^l$& $s\geq 4,\frac{1}{2}s<k<l<\frac{2}{3}s,p\mid k,s,p \nmid l,\lambda \neq 0$\\
&& and do not exist $3\leq u\leq l-1,4\leq v \leq s-1$\\
&&such that any of the condition \eqref{condition-x3} holds \\
\hline
$E_{2t,3t,0}$&$x^3+xy^{2t}+\lambda y^{3t}$& $p \neq 31,t\geq 2,\lambda \neq 0$\\
&& and do not exist $3\leq u\leq 2t-1,4\leq v \leq 4t-2$\\
&&(resp. $4\leq v \leq 4t-1$ if additionally $p \mid t$)\\
&&such that any of the condition \eqref{condition-x3} holds \\
\hline
$E_{2t,3t,l}$&$x^3+xy^{2t}+\lambda y^{3t}+y^l$& $t\geq 2,l>3t,p \nmid l-3t,\lambda \neq 0$\\
&& and do not exist $3\leq u\leq 2t-1,4\leq v \leq l-1$\\
&&such that any of the condition \eqref{condition-x3} holds \\
\hline
$W_{12}$&$x^4+y^5$ &$p \neq 5$ \\
\hline
$W_{12}'$&$x^4+y^5+x^2y^3$&$p \neq 5$ \\
\hline
$W_{13}$&$x^4+xy^4$& $p \neq 5$\\
\hline
$W_{13}'$&$x^4+xy^4+y^6$ & $p \neq 5$\\
\hline
$W_{1,0}$&$x^4+x^2y^3+\lambda y^6$& $\lambda \neq 0,\frac{1}{4},\ p \neq 5$  \\
\hline
$W_{1,0}'$&$x^4+x^2y^3+\lambda y^6+y^7$& $\lambda \neq 0,\frac{1}{4},\ p\neq 5$ \\
\hline
$W_{1,t}$&$x^4+x^2y^3+y^t$& $t \geq 7,\ p\neq 5$\\
\hline
$W_{1,0}^\#$&$x^4+y^6$&$p \neq 5$\\
\hline
$W_{1,0}^{\#'}$&$x^4+x^2y^4+y^6$&$p \neq 5$\\
\hline
$W_{17}$&$x^4+xy^5$&$p \neq 5$\\
\hline
$W_{17}'$&$x^4+xy^5+y^7$&$p \neq 5$\\
\hline
$W_{17}''$&$x^4+xy^5+y^8$&$p \neq 5$\\
\hline
$W_{18}$&$x^4+y^7$&$p \neq 5,7$\\
\hline
$W_{18}'$&$x^4+y^7+x^2y^4$&$p \neq 5,7$\\
\hline
$W_{18}$&$x^4+y^7+x^2y^5$&$p \neq 5,7$\\
\hline
$Z_{0,s}$&$x^3y+y^{s}$&$s \geq 5$\\
&&and do not exist $3\leq u\leq s-1$,$3\leq v\leq s-1$\\
&&such that any of the condition \eqref{condition-x3} and \eqref{condition-x3y} holds\\
\hline
$Z_{r,0}$&$x^3y+xy^{r}$&$r \geq 4$\\
&&and do not exist $3\leq u\leq r-1$,$3\leq v\leq 2r-2$\\
&&such that any of the condition \eqref{condition-x3} and \eqref{condition-x3y} holds\\
\hline
$Z_{r,s}^0$&$x^3y+xy^{r}+y^s$&$s \geq 5, \frac{2s+1}{3}<r\leq s-2$\\
&&(resp.$\frac{2s+1}{3}<r\leq s-1$ if additionally $p \mid s$)\\
&&and do not exist $3\leq u\leq r-1$,$3\leq v\leq s-1$\\
&&(resp. $3\leq u\leq s-1$ if additionally $p \mid 3r-2s-1$)\\
&&such that any of the condition \eqref{condition-x3} and \eqref{condition-x3y} holds\\\hline
$Z_{r,s}^1$&$x^3y+xy^{r}+y^s$&$r \geq 4, 3r-1<2s< 4r,p \nmid 3r-2s-1 $\\
&&and do not exist $3\leq u\leq r-1$,$3\leq v\leq s-1$\\
&&(resp. $3\leq u\leq 3r-s+1$ if additionally $p \mid 3r-2s-1$)\\
&&such that any of the condition \eqref{condition-x3} and \eqref{condition-x3y} holds\\\hline
$Z_{k,s,l}^0$&$x^3y+y^s+\lambda xy^k+xy^l$&$s\geq 5,\frac{2s+1}{3}<k<l\leq s-2,$\\
&&$p\mid 3k-2s-1,p \nmid 3l-2s-1,\lambda \neq 0$\\
&&(resp.$\frac{2s+1}{3}<k<l\leq s-1$ if additionally $p \mid s$)\\
&&and do not exist $3\leq u\leq l-1$,$3\leq v\leq s-1$\\
&&such that any of the condition \eqref{condition-x3} and \eqref{condition-x3y} holds\\
\hline
$Z_{k,s,l}^1$&$x^3y+y^s+\lambda xy^k+xy^l$&$s\geq 5,\frac{1}{2}s<k<l< \frac{2s+1}{3},$\\
&&$p\mid 3k-2s-1,p \nmid 3l-2s-1,\lambda \neq 0$\\
&&and do not exist $3\leq u\leq l-1$,$3\leq v\leq s-1$\\
&&such that any of the condition \eqref{condition-x3} and \eqref{condition-x3y} holds\\
\hline
$Z_{2t,3t,0}$&$x^3y+xy^{2t+1}+\lambda y^{3t+1}$& $p \neq 31,t\geq 2,\lambda \neq 0$\\
&&and do not exist $3\leq u\leq 2t$,$3\leq v\leq 4t$\\
&&such that any of the condition \eqref{condition-x3} and \eqref{condition-x3y} holds\\
\hline
$Z_{2t,3t,l}$&$x^3y+xy^{2t+1}+\lambda y^{3t+1}+y^l$& $t\geq 2,l>3t+1,p \nmid l-3t-1,\lambda \neq 0$\\
&&and do not exist $3\leq u\leq 2t$,$3\leq v\leq l-1$\\
&&such that any of the condition \eqref{condition-x3} and \eqref{condition-x3y} holds\\
\hline
$T_{4,s,2}$&$x^4+x^2y^2+y^s$&$s \geq 5$\\
\hline
$T_{r,s,2}$&$x^r+x^2y^2+y^s$&$r,s \geq 5$\\
\hline
$T_{4,4,2}$&$x^4+\lambda x^2y^2+y^4$&$\lambda^2 \neq 4$\\
\hline
\end{longtable}
    where the condition \eqref{condition-x3} is
    \begin{equation}\label{condition-x3}
    \begin{aligned}
        &A:\ \left\lfloor {\frac{2}{3}v} \right\rfloor+1 \leq u \leq u+p \leq v-3,\ p\mid 3u-2v;\\
        &B:\ p \mid u, p \mid v, 3u<2v<4u, u<u+p\leq 3u-v-2;\\
        &C:\ p \neq 31, u \mathrm{\ is\ even},\frac{3}{2}u+1\leq v\leq 2u-3,p \mid v-\frac{3}{2}u;\\
        &D:\ p=31,u \mathrm{\ is\ even},\frac{3}{2}u+1\leq v,p \mid v-\frac{3}{2}u;
    \end{aligned}
\end{equation}
and the condition \eqref{condition-x3y} is
\begin{equation}\label{condition-x3y}
    \begin{aligned}
        &u\geq 4,v \geq 5 \mathrm{\ and}\\
        &A:\ \left\lfloor {\frac{2v+1}{3}} \right\rfloor+1 \leq u \leq u+p \leq v-3,\ p\mid 3u-2v-1;\\
        &B:\ p \mid 3u-2v-1, u<u+p\leq 3u-v;\\
        &C:\ p \neq 31, u \mathrm{\ is\ odd},\frac{3u+1}{2}u\leq v\leq 2u-3,p \mid v-\frac{3u-3}{2};\\
        &D:\ p=31,u \mathrm{\ is\ odd},\frac{3u+1}{2}u\leq v,p \mid v-\frac{3u-3}{2};
    \end{aligned}
\end{equation}

\renewcommand\arraystretch{1.5}
\begin{longtable}{|c|c|c|}
\caption{}\label{final-table-xyz}
\\
\hline
Symbol&Form& condition  \\
\hline
$T_{3,3,3}$&$x^3+y^3+z^3+\lambda xyz $&$\lambda^3+27\neq 0$ \\
\hline
$T_{r,s,t}$&$x^r+y^s+z^t+xyz$&$\max \{r,s,t\}\geq 4$ \\
\hline
$Q_{0,s}$&$x^3+yz^2+y^{s}$&$s \geq 4$ and do not exist $3\leq u<v\leq s-1$\\
&&(resp. $3\leq u<v\leq s$ if additionally $p \mid s$) \\
&&such that any of the condition \eqref{condition-x3} holds \\
\hline
$Q_{r,0}$&$x^3+yz^2+xy^{r}$&$r \geq 3$ and do not exist $3\leq u\leq r-1,4\leq v\leq 2r-2$ \\
&& (resp. $4\leq v\leq 2r-1$ if additionally $p \mid r$)\\
&&such that any of the condition \eqref{condition-x3} holds \\
\hline
$Q_{r,s}^0$&$x^3+yz^2+y^s+xy^r$& $s\geq 4,\frac{2}{3}s<r \leq s-2,p\nmid 3r-2s$\\
&& (resp. $\frac{2}{3}s<r \leq s-1$ if additionally $p \mid s$)\\
&& and do not exist $3\leq u\leq r-1,4\leq v \leq s-1$\\
&&such that any of the condition \eqref{condition-x3} holds \\
\hline
$Q_{r,s}^{0'}$&$x^3+yz^2+y^s+xy^r$& $s\geq 4,\frac{2}{3}s<r \leq s-2,p\mid 3r-2s$\\
&except for $x^3+xy^4+y^5$& (resp. $\frac{2}{3}s<r \leq s-1$ if additionally $p \mid s$)\\
&when $p=5$& and do not exist $3\leq u\leq s-2,4\leq v \leq s-1$\\
&(which is simple)& (resp. $3\leq u\leq s-2$ if additionally $p \mid s$)\\
&&such that any of the condition \eqref{condition-x3} holds \\
\hline
$Q_{r,s}^1$&$x^3+yz^2+y^s+xy^r$& $r\geq 3,3r<2s<4r,p\nmid 3r-2s$\\
&& and do not exist $3\leq u\leq r-1,4\leq v \leq s-1$\\
&&such that any of the condition \eqref{condition-x3} holds \\
\hline
$Q_{r,s}^{1'}$&$x^3+yz^2+y^s+xy^r$& $r\geq 3,3r<2s<4r,p\mid 3r-2s$\\
&& and do not exist $3\leq u\leq r-1,4\leq v \leq s$\\
&&(resp. $3\leq u \leq 3r-s-1$ if additionally $p \mid r,s$)\\
&&such that any of the condition \eqref{condition-x3} holds \\
\hline
$Q_{k,s,l}^0$&$x^3+yz^2+y^s+\lambda xy^k+xy^l$& $s\geq 4,\frac{2}{3}s<k<l\leq s-2,p\mid 3k-2s,p \nmid 3l-2s,\lambda \neq 0$\\
&&(resp. $\frac{2}{3}s<k<l\leq s-1$ if additionally $p \mid s$)\\
&& and do not exist $3\leq u\leq l-1,4\leq v \leq s-1$\\
&&such that any of the condition \eqref{condition-x3} holds \\
\hline
$Q_{k,s,l}^1$&$x^3+yz^2+y^s+\lambda xy^k+xy^l$& $s\geq 4,\frac{1}{2}s<k<l<\frac{2}{3}s,p\mid k,s,p \nmid l,\lambda \neq 0$\\
&& and do not exist $3\leq u\leq l-1,4\leq v \leq s-1$\\
&&such that any of the condition \eqref{condition-x3} holds \\
\hline
$Q_{2t,3t,0}$&$x^3+yz^2+xy^{2t}+\lambda y^{3t}$& $p \neq 31,t\geq 2,\lambda \neq 0$\\
&& and do not exist $3\leq u\leq 2t-1,4\leq v \leq 4t-2$\\
&&(resp. $4\leq v \leq 4t-1$ if additionally $p \mid t$)\\
&&such that any of the condition \eqref{condition-x3} holds \\
\hline
$Q_{2t,3t,l}$&$x^3+yz^2+xy^{2t}+\lambda y^{3t}+y^l$& $t\geq 2,l>3t,p \nmid l-3t,\lambda \neq 0$\\
&& and do not exist $3\leq u\leq 2t-1,4\leq v \leq l-1$\\
&&such that any of the condition \eqref{condition-x3} holds \\
\hline
$S_{11}$&$x^2z+yz^2+y^4$ & \\
\hline
$S_{11}'$&$x^2z+yz^2+y^4+\lambda x^2y^2$ & \\
\hline
$S_{12}$&$x^2z+yz^2+xy^3$& \\
\hline
$S_{1,0}$&$x^2z+yz^2+x^2y^2+\lambda y^5$&$\lambda \neq 0$\\
\hline
$S_{1,0}^{1}$&$x^2z+yz^2+x^2y^2+\lambda y^5+y^6$&$\lambda \neq 0$\\
\hline
$S_{1,0}^{2}$&$x^2z+yz^2+x^2y^2+xy^4$&\\
\hline
$S_{1,0}^{3}$&$x^2z+yz^2+y^5$ &$p \neq 5$\\
\hline
$S_{1,0}^{4}$&$x^2z+yz^2+x^2y^3+y^5$ &$p \neq 5$\\
\hline
$S_{1,0,t}$&$x^2z+yz^2+x^2y^2+ y^t$&$6\leq t<s+2$\\
\hline
$S_{1,s,0}$&$x^2z+yz^2+x^2y^2+ xy^s$&$t\geq 2s-2$\\
\hline
$S_{1,s,t}$&$x^2z+yz^2+x^2y^2+ xy^s+\lambda y^{t}$&$s\geq 5,s+2\leq t\leq 2s-3,\lambda \neq 0$\\
\hline
$S_{16}$&$x^2z+yz^2+xy^4$&$p \neq 5$ \\
\hline
$S_{16}'$&$x^2z+yz^2+xy^4+y^6$&$p \neq 5$ \\
\hline
$S_{16}''$&$x^2z+yz^2+xy^4+y^7$&$p \neq 5$ \\
\hline
$S_{17}$&$x^2z+yz^2+y^6$ &$p \neq 5$ \\
\hline
$S_{17}'$&$x^2z+yz^2+y^6+x^2y^3$ &$p \neq 5$ \\
\hline
$S_{17}''$&$x^2z+yz^2+y^6+x^2y^4$ &$p \neq 5$ \\
\hline
$U_{12}$&$x^3+xz^2+y^4$&   \\
\hline
$U_{12}'$&$x^3+xz^2+y^4+x^2y^2$& \\
\hline
$U_{1,0}$&$x^3+xz^2+xy^3+\lambda y^3z$& $\lambda^2 \neq 0,-1$  \\
\hline
$U_{1,0}'$&$x^3+xz^2+xy^3+\lambda y^3z+y^4z$& $\lambda^2 \neq 0,-1$  \\
\hline
$U_{1,t}$&$x^3+xz^2+xy^3+y^tz$& $t \geq 4,p \neq 5$  \\
\hline
$U_{16}$&$x^3+xz^2+y^5$& $p \neq 5$  \\
\hline
$U_{16}'$&$x^3+xz^2+y^5+x^2y^3$& $p \neq 5$  \\
\hline
$U_{*}$&$x^3+xz^2+y^3z$&   \\
\hline
$U_{*}'$&$x^3+xz^2+y^3z+xy^4$&   \\
\hline
\end{longtable}

and $g(x_1,x_2)+x_3^2+\dots+x_n^2$ or $h(x_1,x_2,x_3)+x_4^2+\dots+x_n^2$, where $g(x_1,x_2)$ is one of the forms in Table \ref{final-table-xy} and $h(x_1,x_2,x_3)$ is one of the forms in Table \ref{final-table-xyz}.
\end{theorem}

\newcommand{\etalchar}[1]{$^{#1}$}


\begin{thebibliography}{GPB{\etalchar{+}}08}

\bibitem[Arn76]{arnold-class-C}
V.I. Arnold.
\newblock Local normal forms of functions.
\newblock {\em Invent. Math.}, 35:87--109, 1976.


\bibitem[AVGZ12]{arnold}
V.I. Arnold, A.N. Varchenko, and S.M. Gusein-Zade.
\newblock {\em  Singularities of differentiable maps. Volume 1. }
\newblock Reprint of the 1985 edition. {\em Modern Birkhäuser Classics.} Birkhäuser/Springer, New York, 2012. xii+382 pp. ISBN: 978-0-8176-8339-9. 

\bibitem[BGM10]{2010_finite_deter}
Y.~Boubakri, G.M.~Greuel, and T.~Markwig.
\newblock Invariants of hypersurface singularities in positive characteristic.
\newblock {\em Rev. Mat. Complut.}, 25:61--85, 2010.

\bibitem[BGM11]{normal}
Y.~Boubakri, G.M. Greuel, and T.~Markwig.
\newblock Normal forms of hypersurface singularities in positive characteristic.
\newblock {\em Moscow Math.  J.}, 11:657--683, 2011.

\bibitem[DG83]{complete-unimodal-plane}
A.~Dimca and C.G. Gibson.
\newblock Contact unimodular germs from the plane to the plane.
\newblock {\em Q. J. Math.}, 34(3):281--295, 1983.

\bibitem[GK90]{1990_simple}
G.M. Greuel and H.~Kröning.
\newblock Simple singularities in positive characteristic.
\newblock {\em Math.  Z.}, 203:339--354, 1990.

\bibitem[GN16]{right-simple}
G.M. Greuel and H.D. Nguyen.
\newblock Right simple singularities in positive characteristic.
\newblock {\em J. Reine Angew. Math.}, 712:81--106, 2016.

\bibitem[GP17]{mather-yau}
G.M. Greuel and T.H. Pham.
\newblock Mather-{Y}au theorem in positive characteristic.
\newblock {\em J. Algebraic Geom.}, 26(2):347--355, 2017.

\bibitem[GPB{\etalchar{+}}08]{singular-introduction}
G.M. Greuel, G.~Pfister, O.~Bachmann, C.~Lossen, and H.~Schönemann.
\newblock {\em A singular introduction to commutative algebra}.
\newblock Springer-Verlag, 2008.

\bibitem[Mil17]{Milne_alg_group}
J.S. Milne.
\newblock Algebraic groups: the theory of group schemes of finite type over a field.
\newblock {\em Cambridge Studies in Advanced Mathematics. Cambridge University Press, 2017.}

\bibitem[MYZ26]{ma2025}
H.R. Ma, Stephen S.~T. Yau, and H.Q. Zuo.
\newblock Classification of unimodal isolated complete intersection singularities in positive characteristic.
\newblock {\em Pacific J. Math.}, 340(1):71--114, 2026.

\bibitem[Ngu13]{phdclassification}
H.D. Nguyen.
\newblock Classification of singularities in positive characteristic.
\newblock PhD thesis, 2013.

\bibitem[Ngu17]{right_unimodal}
H.D. Nguyen.
\newblock Right unimodal and bimodal singularities in positive characteristic.
\newblock {\em Int. Math. Res. Not. } 2019:1612--1641, 2017.

\bibitem[PG19]{finitedeter}
T.H. Pham and G.M. Greuel.
\newblock On finite determinacy for matrices of power series.
\newblock {\em Math.  Z.}, 290, 2019.

\bibitem[PPG25]{icissimple}
T.H. Pham, G.~Pfister, and G.M. Greuel.
\newblock Classification of simple 0-dimensional isolated complete intersection singularities.
\newblock {\em Int. Math. Res. Not. IMRN}, (14):rnaf191, 22pp.,  2025.

\bibitem[Ros56]{Rosenlicht}
M.~Rosenlicht.
\newblock Some basic theorems on algebraic groups.
\newblock {\em Amer.  J.  Math.}, 78(2):401--443, 1956.

\bibitem[Wal83]{Wall}
C.T.C. Wall.
\newblock Classification of unimodal isolated singularities of complete intersections.
\newblock In {\em Singularities, {P}art 2 ({A}rcata, {C}alif., 1981)}, volume~40 of {\em Proc. Sympos. Pure Math.}, pages 625--640. Amer. Math. Soc., Providence, RI, 1983.

\bibitem[WR05]{actions-algebraic-group}
F.S. Walter and A.~Rittatore.
\newblock  Actions and invariants of algebraic groups.
\newblock CRC Press, 2005.

\bibitem[GPf21]{Semicontinuity}
G.M. Greuel and  G.~Pfister.
Semicontinuity of singularity invariants in families of formal power series. 
Singularities and their interaction with geometry and low dimensional topology-in honor of András Némethi, 207-245, Trends Math., Birkhäuser/Springer, Cham,  2021. 





\end{thebibliography}
\end{document}